\documentclass[11pt]{article}
\usepackage[latin1]{inputenc}
\usepackage{amsmath}
\usepackage{amsfonts}
\usepackage{amssymb}
\usepackage{graphicx}
\usepackage{amsmath}
\usepackage{amsthm}
\usepackage{enumerate}
\usepackage{verbatim}
\usepackage{hyperref}
\usepackage{bbm}
\usepackage{color}

\numberwithin{equation}{section}

\def \R{\mathbb{R}}

\def \E{\mathbb{E}}
\def \P{\mathbb{P}}

\theoremstyle{plain}
\newtheorem{theorem}{Theorem}[section]

\newtheorem{corollary}[theorem]{Corollary}

\newtheorem{lemma}[theorem]{Lemma}

\newtheorem{remark}[theorem]{Remark}

\begin{document}

\title{Non universality for the variance of the number of real roots\\
of random trigonometric polynomials}
\author{\textsc{Vlad Bally}\thanks{
Universit\'e Paris-Est, LAMA (UMR CNRS, UPEMLV, UPEC), MathRisk INRIA,
F-77454 Marne-la-Vall\'ee, France. Email: \texttt{bally@univ-mlv.fr} }
\smallskip \\
\textsc{Lucia Caramellino}\thanks{
Dipartimento di Matematica, Universit\`a di Roma ``Tor Vergata'', and
INDAM-GNAMPA, Via della Ricerca Scientifica 1, I-00133 Roma, Italy. Email: 
\texttt{caramell@mat.uniroma2.it}.}\smallskip\\
\textsc{Guillaume Poly}\thanks{
IRMAR, Universit\'e de Rennes 1, 263 avenue du G\'en\'eral Leclerc, CS 74205
35042 Rennes, France. Email: \texttt{guillaume.poly@univ-rennes1.fr} }}
\maketitle

\begin{abstract}
In this article, we consider the following family of random trigonometric
polynomials $p_n(t,Y)=\sum_{k=1}^n Y_{k,1} \cos(kt)+Y_{k,2}\sin(kt)$ for a given
sequence of i.i.d. random variables $\{Y_{k,1},Y_{k,2}\}_{k\ge 1}$ which are
centered and standardized. We set $\mathcal{N}([0,\pi],Y)$ the number of
real roots over $[0,\pi]$ and $\mathcal{N}([0,\pi],G)$ the corresponding
quantity when the coefficients follow a standard Gaussian distribution. We
prove under a Doeblin's condition on the distribution of the coefficients
that 
\begin{equation*}
\lim_{n\to\infty}\frac{\text{Var}\left(\mathcal{N}_n([0,\pi],Y)\right)}{n}%
=\lim_{n\to\infty}\frac{\text{Var}\left(\mathcal{N}_n([0,\pi],G)\right)}{n}%
+\frac{1}{30}\left(\mathbb{E}(Y_{1,1}^4)-3\right).
\end{equation*}
The latter establishes that the behavior of the variance is not universal
and depends on the distribution of the underlying coefficients through their kurtosis. Actually, a more
general result is proven in this article, which does not requires that the coefficients
are identically distributed. The proof mixes a recent result regarding
Edgeworth's expansions for distribution norms established in \cite{[BCP]}
with the celebrated Kac-Rice formula.
\end{abstract}

\parindent 0pt

\textbf{Keywords}. Random trigonometric polynomials; Edgeworth expansion for non smooth functions; Kac-Rice formula; small ball estimates.

\medskip

\textbf{Mathematics Subject Classifications (2010)}: 60G50, 60F05.

\section{Introduction}

The study of level sets of random functions is a central topic in
probability theory, furthermore at the crossroad of several other domains of
mathematics and physics. In this framework, universality results refer to
asymptotic properties of these random level sets, holding regardless of the
specific nature of the randomness involved. Establishing such universal
properties for generic zero sets allows one to manage what would be
intricate objects. As such, the literature on this topic is very extended
and we refer to the introduction of \cite{[TV14]} and the references therein
for a more exhaustive overview. ~\newline
~\newline
Among the great variety of models that have been investigated, the most
emblematic one is perhaps the so-called Kac polynomials $P_n(x)=\sum_{k=1}^n
a_k x^k$. Assume first that the coefficients $(a_k)_{1\le k\le n}$ are
chosen independently and according to the same centered and standardized
distribution ($\mathbb{E}(a_1)=0,\,\mathbb{E}(a_1^2)=1$). Then, set $%
\mathcal{N}_n(\mathbb{R})$ its number of real roots:

\begin{equation*}
\mathcal{N}_n(\mathbb{R})=\text{card}\left\{x\in\mathbb{R}%
\,\left|\right.\,P_n(x)=0\right\}.
\end{equation*}

As a synthesis of the following (non exhaustive) list of landmark articles 
\cite{[Kac48],[Erd56],[IM71],[Mas74]} the following phenomena hold under
mild conditions, \textit{universally}, that is to say regardless of the
choice of the peculiar distribution of the coefficients:

\medskip

\noindent $\bullet$ \textbf{universality of the mean}: $\displaystyle\mathbb{%
E}\left(\mathcal{N}_n(\mathbb{R})\right)\sim \frac{2}{\pi}\log(n);$

\smallskip

\noindent $\bullet$ \textbf{universality of the variance}: $\displaystyle 
\text{Var}\left(\mathcal{N}_n(\mathbb{R})\right)\sim\frac{4}{\pi}\left(1-%
\frac{2}{\pi}\right)\log(n);$

\smallskip

\noindent $\bullet$ \textbf{universality of the fluctuations around the mean:%
} $\displaystyle \frac{\mathcal{N}_n(\mathbb{R})-\mathbb{E}\left(\mathcal{N}%
_n(\mathbb{R})\right)}{\sqrt{\text{Var}\left(\mathcal{N}_n(\mathbb{R})\right)%
}}\xrightarrow[n\to\infty]{\text{Law}}\mathcal{N}(0,1).$

\medskip

Above, the notation $u_n\sim v_n$ means $\frac{u_n}{v_n}\to 1$ as $%
n\to\infty $, and $\mathcal{N}(0,1)$ stands for the standard normal law.
Many other models of random polynomials exist in the literature for which
universal properties have been intensively investigated. For most of them,
both local universality (i.e. joint distribution of roots at microscopic
scales) and universality of the expectation at a global scale have been
achieved successfully. Concerning local universality, we refer to \cite%
{[TV14],[DNV], [IKM]} and for expectation to \cite{[NNV16],[Fla16],[FK17]}.
Very often, the extension to the global scale of the microscopic
distribution of the roots is not an easy task, and one needs first to
provide suitable estimates for the so-called phenomenon of repulsion of
zeros. Let us also mention that multivariate models have been recently
studied, for which we refer to \cite{[APV17],[CNNV17]}. To the best of our
knowledge, it must be emphasized that the universality of the variance has
only been reached for Kac polynomials.

\smallskip

Here, we investigate this problem for trigonometric models and show that the
variance behavior is actually \textit{not universal} by computing exactly
the correction with respect to the case of Gaussian coefficients. This
results displays a strong difference with the well-known Kac polynomials
models. We stress that our main result only requires the independence of the
coefficients. More concretely, we shall consider for different sequences of
independent random vectors $Y_k=(Y_{k}^1,Y_k^{2})$, $k\in\mathbb{N}$, the
number $\mathcal{N}(0,\pi)$ of real roots over the set $[0,\pi]$ of 
\begin{equation*}
p_n(t,Y)=\sum_{k=1}^n Y_k^{1} \cos( k t) +Y_k^{2} \sin (k t).
\end{equation*}
In order to take benefit from the Central Limit Theorem (hereafter, CLT), we
first make a scale change and rather consider 
\begin{equation*}
P_n(t,Y)=\frac{1}{\sqrt{n}} \sum_{k=1}^n Y_k^{1} \cos\left(\frac{k t}{n}%
\right)+Y_k^{2} \cos\left(\frac{k t}{n}\right), \quad t\in[0,2n\pi].
\end{equation*}
Indeed, it can be established that $P_n(\cdot,Y)$ converge in distribution
towards a stationary Gaussian process whose correlation function is $\frac{%
\sin (x)}{x}$. On the other hand, doing so, the number of roots of $p_n(\cdot,Y)$ over $[0,\pi]$
is also the number of roots of $P_n(\cdot,Y)$ over $[0,n\pi]$ and one loses
nothing in this procedure. We also highlight that $P_n(\cdot,Y)$ is much
more manageable thanks to the aforementioned limit theorem.

\medskip

In order to state more precisely our main theorem, we need some preliminary
notations given in the following subsection.

\medskip

\textbf{Main result.} We consider a sequence of centered, independent random
vectors $\{Y_{k}\}_{k\geq 1}\in \mathbb{R}^{2}$ with the normalization $\mathbb{E}(Y_{k}^{i}Y_{k}^{j})=\delta _{i,j}$  and which satisfy Doeblin condition (\ref{D}) with the moment conditions (\ref{M}). Next, we consider the following trigonometric polynomials:
\begin{equation}
P_{n}(t,Y)=\frac{1}{\sqrt{n}}\sum_{k=1}^{n}\cos (\frac{kt}{n})Y_{k}^{1}+\sin
(\frac{kt}{n})Y_{k}^{2}  \label{pol}
\end{equation}%
and we denote by $N_{n}(Y)$ the number of roots of $P_{n}(t,Y)$ in the
interval $(0,n\pi )$. We shall focus on the variance of $N_{n}(Y)$ given by 
\begin{equation*}
\text{Var}\left( N_{n}(Y)\right) =\mathbb{E}(N_{n}^{2}(Y))-(\mathbb{E}%
(N_{n}(Y))^{2}
\end{equation*}

It is known thanks to the appearance of \cite{[GW11]} that if $G=(G_{k})_{k\in
N} $ is a sequence of two dimensional standard random variables then the
following limit exists 
\begin{equation*}
\lim_{n}\frac{1}{n}\text{Var}\left( N_{n}(G)\right) =C(G)\approx 0.56
\end{equation*}%
(for the explicit expression of $C(G)$ see page 298 of \cite{[GW11]}, we stress that the previous approximation of $C(G)$ concerns the number of zeros over $[0,2\pi]$).
Besides a Central Limit Theorem is also established regarding the
fluctuations of the number of roots around the mean. We also refer to \cite%
{[AL13],[ADL16]} for alternative proofs and some refinements
obtained by following the so-called \textit{Nourdin-Peccati method} for
establishing central limit theorems for functionals of Gaussian processes.
Our aim is to prove a similar result for the variance of $N_{n}(Y) $ and all the more to
compute explicitly the constant $C(Y).$ At this point, it must be emphasized
that outside the scope of functionals of Gaussian processes, one cannot
anymore deploy the powerful combination of Malliavin calculus and Wiener
chaos theory as explained in the book \cite{[NP12]}. In order to bypass this
restriction, as explained below, our approach heavily relies on combination
of Edgeworth expansion and Kac-Rice formulae. Let us also mention that the
universality of the expected number of roots has been recently fully
established in \cite{[Fla16]} under a second moment condition.

\bigskip

An important aspect of our contribution is that we can formulate explicitly $C(Y)$. Our main result is the following (see Theorem \ref{Main-result}). Suppose that $Y\in 
\mathcal{D}(\varepsilon ,r)$ and suppose also that for every multi-index $%
\alpha $ with $\left\vert \alpha \right\vert =3,4$\ the following limits
exists and are finite: 
\begin{equation*}
\lim_{n}\E(\prod_{i=1}^{3}Y_{n}^{\alpha _{i}})=y_{\infty }(\alpha ),\quad
\lim_{n}\E(\prod_{i=1}^{4}Y_{n}^{\alpha _{i}})=y_{\infty }(\alpha ).
\end{equation*}
Then 
\begin{equation*}
\lim_{n}\frac{1}{n}V_{n}(Y)=C(G)+\frac{1}{60}\times y_{\ast }
\end{equation*}
with 
\begin{equation*}
y_{\ast }=\left( (y_{\infty }(1,1,2,2)-1)+(y_{\infty
}(2,2,1,1)-1)+(y_{\infty }(1,1,1,1)-3)+(y_{\infty }(2,2,2,2)-3)\right) .
\end{equation*}

Notice that the random vectors $(Y_{k})_{k\geq 1}$ are not supposed here to be
identically distributed. However, the hypothesis (\ref{D}) and (\ref{M})
display some uniformity because $\eta ,r$ and $M_{p}(Y)$ are uniform
parameters. For simplicity, suppose for a moment that they are uniformly
distributed and moreover that the components $Y_{k}^1$ and $Y_{k}^2$ of $%
Y_{k}=(Y_{k}^ 1,Y_{k}^2$) are also i.i.d. Then $y_{\infty
}(1,1,2,2)=y_{\infty }(2,2,1,1)=1$ and $y_{\infty }(1,1,1,1)=y_{\infty
}(2,2,2,2)=\mathbb{E}((Y^1_1)^{4})$. In such a case, the non-universality of
the variance becomes more transparent since 
\begin{equation}
\lim_{n\rightarrow \infty }\frac{\text{Var}\left( N_{n}(Y)\right) }{n}%
=\lim_{n\rightarrow \infty }\frac{\text{Var}\left( N_{n}(G)\right) }{n}%
+\frac{1}{30}\left( \mathbb{E}\left( (Y_1^1)^{4}\right) -3\right).
\label{univareqiid}
\end{equation}%
In particular, one notice that the deviation from the Gaussian behavior is
exactly proportional to the \textit{kurtosis} of the random variables under
consideration.

\medskip

\textbf{Strategy of the proof.} Let us summarize briefly the main steps of
our proofs. Basically, up to some technical details, it illustrates rather
well the main ideas of our approach.

\bigskip

\underline{Step 1: An approximated Kac-Rice formula}

\bigskip

Let us recall the celebrated and very useful Kac-Rice formula. Consider a
smooth deterministic function $f$ defined on $[a,b]$ such that $%
|f(t)|+|f^{\prime }(t)|>0$ for all $t\in \lbrack a,b]$. Then, one has

\begin{equation*}
\text{Card}\left\{t\in]a,b[\,\,\left|\right.\,\,f(t)=0\right\}=\lim_{\delta%
\to 0} \frac{1}{2\delta} \int_a^b |f^{\prime }(t)|\mathbf{1}%
_{\left\{|f(t)|<\delta\right\}}dt.
\end{equation*}

When one applies the latter to the random functions $P_n(t,Y)$ one needs to
handle the level of \textit{non degeneracy} which determines the speed of
convergence in the Kac-Rice formula. More concretely, in our proof, we will
use that for $\delta_n=\frac{1}{n^5}$:

\begin{equation*}
\lim_{n\rightarrow \infty }\frac{1}{n}\text{Var}\left(N_n(Y)\right)
-\lim_{n\rightarrow \infty }\frac{1}{n}\text{Var}\left( \frac{1}{2\delta _{n}%
}\int_{0}^{n\pi }|P_{n}^{\prime }(t,Y)|\mathbf{1}_{\left\{
|P_{n}(t,Y)|<\delta _{n}\right\} }dt\right) =0.
\end{equation*}

We refer to Lemma \ref{KR} for this step.

\bigskip

\underline{Step 2: Removing the diagonal}

\bigskip

When computing the variance, expressions of the following kind appear:%
\begin{equation}  \label{kr}
\begin{array}{l}
\displaystyle \frac{1}{n}\int_{0}^{n\pi }\int_{0}^{n\pi }\Phi
_{n}(t,s,Y)dsdt, \mbox{ where}\smallskip \\ 
\displaystyle \Phi _{n}(t,s,Y) =|P_{n}^{\prime }(t,Y)|\frac{1}{2\delta _{n}}%
\mathbf{1}_{\left\{ |P_{n}(t,Y)|<\delta _{n}\right\} }|P_{n}^{\prime }(s,Y)|%
\frac{1}{ 2\delta _{n}}\mathbf{1}_{\left\{ |P_{n}(s,Y)|<\delta _{n}\right\}
} .%
\end{array}%
\end{equation}
Notice that 
\begin{equation*}
P_{n}^{\prime }(t,Y)=\frac{1}{\sqrt{n}}\sum_{k=1}^{n}\frac{k}{n}\cos (\frac{%
kt}{n})Y_{k}^{2}-\frac{k}{n}\sin (\frac{kt}{n})Y_{k}^{1}
\end{equation*}%
so it becomes clear that in order to study the asymptotic behaviour of $%
E(\Phi _{n}(t,s,Y))$ one has to use the CLT for the random vector $%
S_{n}(t,s,Y)=(P_{n}(t,Y),P_{n}^{\prime }(t,Y),P_{n}(s,Y),P_{n}^{\prime
}(s,Y)).$ A first difficulty in doing this is that $\frac{1}{2\delta _{n}}%
\mathbf{1}_{\left\{ |P_{n}(t,Y)|<\delta _{n}\right\} }\rightarrow \delta
_{0}(P_{n}(t,Y))$ so we are out from the framework of continuous and bounded
test functions considered in the classical CLT. We have to use a variant of
this theorem concerning convergence in distribution norms -- this result is
established in \cite{[BCP]}. A second difficulty concerns the non degeneracy
of the vector $S_{n}(t,s,Y)$: when $t\approx s$, the random vector $%
S_{n}(t,s,Y)$ becomes degenerate and employing the central limit theorem or
its Edgeworth's expansions turn out to be hard. In order to avoid that, we
give us a fixed parameter $\epsilon >0$ and we prove that
\begin{equation*}
\lim_{\epsilon\to 0} \limsup_{n\to\infty} \frac{1}{n}\left( \frac{1}{4
\delta_n^2} \int_{[0,n\pi]^2, |t-s|< \epsilon} \text{Cov}%
\left(|P_n^{\prime }(t,Y)|\mathbf{1}_{\left\{|P_n(t,Y)|<\delta_n\right%
\}},|P_n^{\prime }(s,Y)|\mathbf{1}_{\left\{|P_n(s,Y)|<\delta_n\right\}}%
\right)dt ds\right) =0.
\end{equation*}
The latter enables us to impose the condition $|t-s|\geq \epsilon $ in all
our Kac-Rice estimates. This is particularly convenient since the underlying
processes become uniformly non-degenerate.

A third difficulty comes for the fact that, roughly speaking,%
\begin{equation*}
\frac{1}{n}\int_{0}^{n\pi }\int_{0}^{n\pi }\big | \E(\Phi _{n}(t,s,Y))-\E(\Phi
_{n}(t,s,G))\big | dsdt\sim \frac{1}{n}\times (\pi n)^{2}\times \big | %
\E(\Phi _{n}(\cdot,\cdot,Y))-\E(\Phi _{n}(\cdot ,\cdot,G))\big |
\end{equation*}%
so we need to get 
\begin{equation*}
\left\vert \E(\Phi _{n}(t,s,Y))-\E(\Phi _{n}(t,s,G))\right\vert \leq \frac{C}{%
n^{3/2}}
\end{equation*}%
and in order to achieve this, it is not sufficient to use the CLT, but we
have to use an Edgeworth expansion of order three.

\bigskip

\underline{Step 3: Performing Edgeworth's expansions}

\bigskip

In this step, we make use of Edgeworth's expansion in distribution norm
developed in \cite{[BCP]}. We first set

\begin{equation*}
F_n(x_1,x_2,x_3,x_4)=\frac{1}{4\delta_n^2}\times |x_1|\mathbf{1}%
_{\{|x_2|<\delta_n}|x_3|\mathbf{1}_{\{|x_4|<\delta_n},
\end{equation*}

and $\rho_{n,t,s} $ the density of $(P_n(t,G),P_n^{\prime
}(t,G),P_n(s,G),P_n^{\prime }(s,G))$. By using the Edgeworth's expansion, we
will prove that

\begin{eqnarray*}
&&\mathbb{E}\left( F_n\left(P_n(t,Y),P_n^{\prime }(t,Y),P_n(s,Y),P_n^{\prime
}(s,Y)\right)\right)= \\
&&\int_{\mathbb{R}^4} F_n(x) \rho_{n,t,s}(x)\left(1+\frac{1}{\sqrt{n}}
Q_{n,t,s}(x)+\frac{1}{n} R_{n,t,s}(x) \right) dx_1 dx_2 dx_3 dx_4 \\
&& + \mathcal{R}_n(t,s).
\end{eqnarray*}

where $Q_n$ and $R_n$ are totally explicit polynomials of degree less than $%
6 $ whose coefficients involve the moments of the sequence of the random
variables $\{Y_k^1,Y_k^2\}_{k\ge 1}$ and where the remaining term satisfies

\begin{equation*}
\lim_{\epsilon\to 0}\lim_{n\to\infty} \frac{1}{n}\int_{[0,n\pi]^2,|t-s|\ge
\epsilon} \mathcal{R}_n(t,s) dt ds =0.
\end{equation*}

Doing so, some computations are involved but they are totally transparent in
terms of the moments of the coefficients of our polynomial. This step allows
one to handle explicitly the various cancellations occurring in the
variance. This step is the heart of the proof and is done in Section \ref%
{Cancellations}. We strongly emphasize that getting a polynomial speed of
convergence in the Kac-Rice formula is crucial in order to manage the
remainder of the Edgeworth's expansions.

\section{The problem}

We consider a sequence of centered, independent random variables $Y_{k}\in
\R^{2},k\in N$ with $\E(Y_{k}^{i}Y_{k}^{j})=\delta _{i,j}$. We assume that
they satisfy the following ``Doeblin condition'': there exists some points $%
y_{k}\in \R^{2}$ and $r,\eta \in (0,1)$ such that for every $k\in N$ and
every measurable set $A\subset B_{r}(y_{k})$ 
\begin{equation}
\P(Y_{k}\in A)\geq \eta \mathrm{Leb}(A).  \label{D}
\end{equation}%
Moreover we assume that $Y_{k},k\in N$ have finite moments of any order
which are uniformly bounded with respect to $k:$%
\begin{equation}
\sup_{k}(\E(\left\vert Y_{k}\right\vert ^{p}))^{1/p}=M_{p}(Y)<\infty .
\label{M}
\end{equation}%
We denote by $\mathcal{D}(r,\eta )$ the sequences of random variables $%
Y=(Y_{k})_{k\in N}$ which are independent and verify (\ref{D}) and (\ref{M})
for every$p\geq 1$. Moreover we put 
\begin{equation}
P_{n}(t,Y)=\frac{1}{\sqrt{n}}\sum_{k=1}^{n}\cos (\frac{kt}{n})Y_{k}^{1}+\sin
(\frac{kt}{n})Y_{k}^{2}  \label{pol}
\end{equation}%
and we denote by $N_{n}(Y)$ the number of roots of $P_{n}(t,Y)$ in the
interval $(0,n\pi )$ and by $V_{n}(Y)$ the variance of $N_{n}(Y):$%
\begin{equation}
V_{n}(Y)=\E(N_{n}^{2}(Y))-(\E(N_{n}(Y))^{2}  \label{var}
\end{equation}%
It is known (see e.g. \cite{[GW11]}) that if $G=(G_{k})_{k\in N}$ is a
sequence of two dimensional standard random variables then the following
limit exists 
\begin{equation*}
\lim_{n}\frac{1}{n}V_{n}(G)=C(G).
\end{equation*}%
Our main result is the following.
\begin{theorem}\label{Main-result}
Suppose that $Y\in \mathcal{D}(\eta ,r)$ and suppose also that for every
multi-index $\alpha $ with $\left\vert \alpha \right\vert =3,4$\ the
following limits exists and are finite:%
\begin{equation*}
\lim_{n}\E(\prod_{i=1}^{3}Y_{n}^{\alpha _{i}})=y_{\infty }(\alpha ),\quad
\lim_{n}\E(\prod_{i=1}^{4}Y_{n}^{\alpha _{i}})=y_{\infty }(\alpha ).
\end{equation*}%
Then 
\begin{equation*}
\lim_{n}\frac{1}{n}V_{n}(Y)=C(G)+\frac{1}{60}\times y_{\ast }
\end{equation*}%
with 
\begin{equation*}
y_{\ast }=\left( (y_{\infty }(1,1,2,2)-1)+(y_{\infty
}(2,2,1,1)-1)+(y_{\infty }(1,1,1,1)-3)+(y_{\infty }(2,2,2,2)-3)\right) .
\end{equation*}
\end{theorem}

\textbf{Proof.} The proof is an immediate consequence of Lemma \ref{KR}
point C (see (\ref{d9})) and of Lemma \ref{Cancel} (see (\ref{N4'})).

\bigskip

\begin{remark}
Notice that the random variables $Y_{k}\in R^{2},k\in N$ are not supposed to
be identically distributed. However, the hypothesis (\ref{D}) and (\ref{M})
contain some uniformity assumptions because $\varepsilon ,r$ and $M_{p}(Y)$
are common for all of them. Suppose for a moment that they are identically
distributed and moreover, that the components $Y^{1}=Y_{k}^{1}$ and $%
Y^{2}=Y_{k}^{2}$ are independent. Then $y_{\infty }(1,1,2,2)=y_{\infty
}(2,2,1,1)=1$ and $y_{\infty }(1,1,1,1)=\E(\left\vert Y^{1}\right\vert ^{4})$
and $y_{\infty }(2,2,2,2)=\E(\left\vert Y^{2}\right\vert ^{4}),$ so $y_{\ast
} $ is the sum of the kurtosis of $Y^{1}$ and of $Y^{2}.$ Put it otherwise:
take $\overline{Y}=((Y^{1})^{2},(Y^{2})^{2}).$ Then $y_{\ast }=0$ iff the
covariance matrix of $\overline{Y}$ coincides with the covariance matrix of
the corresponding $\overline{G}.$
\end{remark}

\bigskip

\section{CLT and Edgeworth expansion}

The main tool in the this paper is the $CLT$ and the Edgeworth development
of order two that we proved in \cite{[BCP]} Proposition 2.5. We recall them
here. We consider a sequence of matrices $C_{n}(k)\in \mathcal{M}_{d\times
2},n,k\in N $ which verify%
\begin{equation}
\Sigma _{n}:=\frac{1}{n}\sum_{k=1}^{n}C_{n}(k)C_{n}^{\ast }(k)\geq
\varepsilon _{\ast }\quad and\quad \sup_{n,k\in N}\left\Vert
C_{n}(k)\right\Vert <\infty .  \label{Nd}
\end{equation}%
We denote 
\begin{equation*}
X_{n,k}=C_{n}(k)Y_{k},\quad G_{n,k}=C_{n}(k)G_{k}
\end{equation*}%
where $Y=(Y_{k})_{k\in N}$ is the sequence introduced in the previous
section and $G=(G_{k})$ is a sequence of independent standard normal random
variables in $R^{2}.$ For a multi-index $\alpha =(\alpha _{1},...,\alpha
_{m})\in \{1,...,d\}^{m}$ we denote $\left\vert \alpha \right\vert =m$ and 
\begin{eqnarray}
\Delta _{\alpha }(X_{n,k}) &=&\E(X_{n,k}^{\alpha })-\E(G_{n,k}^{\alpha
})=\E(\prod_{i=1}^{m}X_{n,k}^{\alpha _{i}})-\E(\prod_{i=1}^{m}G_{n,k}^{\alpha
_{i}})\quad and  \label{N5} \\
c_{n}(\alpha ,X) &=&\frac{1}{n}\sum_{k=1}^{n}\Delta _{\alpha }(X_{n,k})
\label{N5a}
\end{eqnarray}%
By hypothesis, for $\left\vert \alpha \right\vert =1,2$ we have $\Delta
_{\alpha }(X_{n,k})=0.$

For a function $f\in C_{pol}^{\infty }(R^{d})$ ($C^{\infty }$ functions with
polynomial growth) and for $q\in N$ we define $L_{q}(f)$ and $l_{q}(f)$ to
be two numbers such that 
\begin{equation}
\sum_{\left\vert \gamma \right\vert \leq q}\left\vert \partial ^{\gamma
}f(x)\right\vert \leq L_{q}(f)(1+\left\vert x\right\vert )^{l_{q}(f)}.
\label{N5b}
\end{equation}%
Moreover we denote%
\begin{eqnarray*}
S_{n}(Y) &=&\frac{1}{\sqrt{n}}\sum_{k=1}^{n}X_{n,k}=\frac{1}{\sqrt{n}}%
\sum_{k=1}^{n}C_{n}(k)Y_{k}\quad and \\
S_{n}(G) &=&\frac{1}{\sqrt{n}}\sum_{k=1}^{n}G_{n,k}=\frac{1}{\sqrt{n}}%
\sum_{k=1}^{n}C_{n}(k)G_{k}.
\end{eqnarray*}%
The $CLT$ says that, if $Y\in \mathcal{D}(\eta ,r)$ then, for every
multi-index $\gamma $ with $\left\vert \gamma \right\vert \leq q$ 
\begin{eqnarray}
\E(\partial ^{\gamma }f(S_{n}(Y))) &=&\E(\partial ^{\gamma }f(S_{n}(G)))+\frac{%
1}{n^{1/2}}R_{n}(f)\quad with  \label{CLT} \\
\left\vert R_{n}(f)\right\vert &\leq &C(L_{0}(f)+n^{1/2}L_{q}(f)e^{-cn})
\label{CLT1}
\end{eqnarray}%
where $C\geq 1\geq c>0$ are constants which depend on $\eta ,r$ in (\ref{D}%
), on $\varepsilon _{\ast }$ in (\ref{Nd}) and on $M_{p}(Y)$ for a
sufficiently large $p.$

We go further and we recall the Edgeworth development. We consider the
Hermite polynomials $H_{\alpha }$ which are characterized by the equality%
\begin{equation}
\E(\partial ^{\alpha }f(W))=\E(f(W)H_{\alpha }(W))\quad \forall f\in
C_{pol}^{\infty }(R^{d})  \label{H1}
\end{equation}%
where $W\in R^{d}$ is a standard normal random variable. Let us mention that 
$H_{\alpha }$ may be represented as follows. Let $h_{k}$ be the Hermite
polynomial of order $k$ on $R$\ - see Nualart [N] for example, for the
definition and the recurrent construction of $h_{k}.$ Now, for the
multi-index $\alpha $ and for $j\in \{1,...,d\}$ we denote $i_{j}(\alpha
)=card\{i:\alpha _{i}=j\}.$ Then $H_{\alpha
}(x_{1},...,x_{d})=h_{i_{1}(\alpha )}(x_{1})\times ...\times h_{i_{d}(\alpha
)}(x_{d}).$ It is known that $h_{k}$ is even (respectively odd) if $k$ is
even (respectively odd) so $H_{\alpha }$ itself has the corresponding
properties on each variable (we will use this in the sequel).

We introduce now the following functions which represent the correctors of
order one and two in the Edgeworth development:%
\begin{eqnarray}
\Gamma _{n,1}(X,x) &=&\frac{1}{6}\sum_{\left\vert \beta \right\vert
=3}c_{n}(\beta ,X)H_{\beta }(x)\quad and  \label{N7b} \\
\Gamma _{n,2}(X,x) &=&\Gamma _{n,2}^{\prime }(X,x)+\Gamma _{n,2}^{\prime
\prime }(X,x)  \label{N7}
\end{eqnarray}%
and%
\begin{eqnarray}
\Gamma _{n,2}^{\prime }(X,x) &=&\frac{1}{24}\sum_{\left\vert \beta
\right\vert =4}c_{n}(\beta ,X)H_{\beta }(x),\quad  \label{N7a} \\
\Gamma _{n,2}^{\prime \prime }(X,x) &=&\frac{1}{72}\sum_{\left\vert \rho
\right\vert =3}\sum_{\left\vert \beta \right\vert =3}c_{n}(\beta
,X)c_{n}(\rho ,X)H_{(\beta ,\rho )}(x)  \label{N7aa}
\end{eqnarray}%
In Proposition 2.6 from \cite{[BCP]} we prove the following: Let $N\in N.$
For every $f\in C_{pol}^{\infty }(R^{d})$ and for every multi-index $\gamma $
with $\left\vert \gamma \right\vert \leq q$%
\begin{eqnarray}
\E(\partial ^{\gamma }f(S_{n}(Y))) &=&\E(\partial ^{\gamma }f(\Sigma
_{n}^{1/2}W)(1+\frac{1}{\sqrt{n}}\Gamma _{n,1}(\Sigma _{n}^{-1/2}X,W))+\frac{%
1}{n}\Gamma _{n,2}(\Sigma _{n}^{-1/2}X,W))  \label{6'} \\
&&+\frac{1}{n^{3/2}}R_{n}(f).  \notag
\end{eqnarray}%
Here $W\in R^{d}$ is a standard normal random variable and $\Sigma
_{n}^{-1/2}X=(\Sigma _{n}^{-1/2}X_{k})_{k\in N}.$ The remainder $R_{n}(f)$
verifies%
\begin{equation}
\left\vert R_{n}(f)\right\vert \leq C(L_{0}(f)+n^{3/2}L_{q}(f)e^{-cn})
\label{N7c}
\end{equation}%
where $C\geq 1\geq c>0$ are constants which depend on $\eta ,r$ in (\ref{D}%
), on $\varepsilon _{\ast }$ in (\ref{Nd}) and on $M_{p}(Y)$ for a
sufficiently large $p.$

Let us mention some more facts which will be useful in our framework. We
will work with an even function $f$ (so $f(x)=f(-x)).$ Since $W$ and $-W$\
have the same law, and the Hermite polynomials of order three are odd we have%
\begin{equation*}
\E(\partial ^{\gamma }f(\Sigma _{n}^{1/2}W)\Gamma _{n,1}(\Sigma
_{n}^{-1/2}X,W)))=0
\end{equation*}%
so this term does no more appear in our development. Moreover consider a
diagonal matrix $I_{d}(\lambda )$ such that $I_{d}^{i}(\lambda )=\lambda
_{i} $\ and such that $\lambda _{i}\geq \varepsilon _{\ast }.$ Then a
straightforward computation (using the non degeneracy of $\Sigma _{n}$ and
of $I_{d}(\lambda )$ and some standard integration by parts techniques)
gives 
\begin{equation*}
\E(\partial ^{\gamma }f(\Sigma _{n}^{1/2}W)\Gamma _{n,2}(\Sigma
_{n}^{-1/2}X,W)))=\E(\partial ^{\gamma }f(I_{d}^{1/2}(\lambda )W)\Gamma
_{n,2}(I_{d}^{-1/2}(\lambda )X,W)))+r_{n}(f)
\end{equation*}%
with 
\begin{equation}
\left\vert r_{n}(f)\right\vert \leq CL_{0}(f)\times \left\Vert \Sigma
_{n}-I_{d}(\lambda )\right\Vert  \label{N7d}
\end{equation}%
with $C$ depending on $\varepsilon _{\ast }.$ Noticing that the law of $%
S_{n}(G)$ coincides with $\Sigma _{n}^{1/2}W$ we write (\ref{6'}) as%
\begin{eqnarray}
\E(\partial ^{\gamma }f(S_{n}(Y))) &=&\E(\partial ^{\gamma }f(S_{n}(G)))+\frac{%
1}{n}\E(\partial ^{\gamma }f(I_{d}^{1/2}(\lambda )W)\Gamma
_{n,2}(I_{d}^{-1/2}(\lambda )X,W)))  \label{N7e} \\
&&+\frac{1}{n}r_{n}(f)+\frac{1}{n^{3/2}}R_{n}(f).  \notag
\end{eqnarray}%
This is the equality that we will use in the sequel.

We finish this section by recalling a result concerning small balls obtained
in \cite{[BCP]} for $S_{n}(t,Y)=(P_{n}(t,Y),P_{n}^{\prime }(t,Y)).$ So we
have to check the hypothesis in Theorem 3.2 in \cite{[BCP]}. In our case we
have $m=1,d=2$ and 
\begin{equation*}
C_{n}(k,t)=\left( 
\begin{tabular}{ll}
$\cos (\frac{kt}{n})$ & $\sin (\frac{kt}{n})$ \\ 
$-\frac{k}{n}\sin (\frac{kt}{n})$ & $\frac{k}{n}\cos (\frac{kt}{n})$%
\end{tabular}%
\right) .
\end{equation*}%
For every $\xi \in R^{2}$ one has $\left\vert C_{n,k}^{\ast }(t)\xi
\right\vert ^{2}=\xi _{1}^{2}+\frac{k^{2}}{n^{2}}\xi _{2}^{2}\geq \frac{k^{2}%
}{n^{2}}\left\vert \xi \right\vert ^{2}$ so that 
\begin{equation*}
\frac{1}{n}\sum_{k=1}^{n}\left\vert C_{n}(k,t)\xi \right\vert ^{2}\geq \frac{%
1}{n}\sum_{k=1}^{n}\frac{k^{2}}{n^{2}}\left\vert \xi \right\vert ^{2}\geq
\int_{0}^{1}x^{2}dx\times \left\vert \xi \right\vert ^{2}=\frac{1}{3}%
\left\vert \xi \right\vert ^{2}.
\end{equation*}%
This means than the hypothesis (3.9)  in \cite{[BCP]} holds with $\lambda
_{\ast }=\frac{1}{3}$ and we are able to use (3.11) (with $l=a=1,d=2)$ from 
\cite{[BCP]}. Then, for every $\theta >1$ and $\varepsilon >0$ we have 
\begin{equation}
\P(\inf_{t\leq n}\left\vert S_{n}(t,Y)\right\vert \leq n^{-\theta })\leq 
\frac{C}{n^{\theta (d-l)-al-\varepsilon }}=\frac{C}{n^{\theta -1-\varepsilon
}}.  \label{SB}
\end{equation}%
Moreover, by (3.10) in the same theorem in \cite{[BCP]} we get%
\begin{equation}
\sup_{t\geq 0}\P(\left\vert S_{n}(t,Y)\right\vert \leq n^{-\theta })\leq C(%
\frac{1}{n^{2\theta }}+e^{-cn}).  \label{SBB}
\end{equation}

\section{Estimates based on Kac-Rice formula}

In this section we will use Kac-Rice lemma that we recall now. Let $%
f:[a,b]\rightarrow R$ be a differentiable function and let 
\begin{equation}
\omega _{a,b}(f)=\inf_{x\in \lbrack a,b]}(\left\vert f(x)\right\vert
+\left\vert f^{\prime }(x)\right\vert )\quad and\quad \delta _{a,b}(f)=\min
\{\left\vert f(a)\right\vert ,\left\vert f(b)\right\vert ,\omega _{a,b}(f)\}.
\label{d3c}
\end{equation}%
We denote by $N_{a,b}(f)$ the number of solutions of $f(t)=0$ for $t\in
\lbrack a,b]$. The Kac-Rice lemma says that if $\delta _{a,b}(f)>0$ then 
\begin{equation}
N_{a,b}(f)=I_{a,b}(\delta ,f):=\int_{a}^{b}\left\vert f^{\prime
}(t)\right\vert 1_{\{\left\vert f(t)\right\vert \leq \delta \}}\frac{dt}{%
2\delta }\quad for\quad 0<\delta \leq \delta _{a,b}(f).  \label{d3a}
\end{equation}%
Notice that we also have, for every $\delta >0$%
\begin{equation}
I_{a,b}(\delta ,f)\leq 1+N_{a,b}(f^{\prime })  \label{d3b}
\end{equation}%
Indeed, we may assume that $N_{a,b}(f^{\prime })=p<\infty $ and then we take 
$a=a_{0}\leq a_{1}<....<a_{p}\leq a_{p+1}=b$ to be the roots of $f^{\prime
}. $ Since $f$ is monotonic on each $(a_{i},a_{i+1})$ one has $%
I_{a_{i},a_{i+1}}(\delta ,f)\leq 1$ so (\ref{d3b}) holds. In the following
we will refer this result as the K-R lemma.

We will use this formula for $f(t)=P_{n}(t,Y).$ We denote%
\begin{equation}
\phi _{\delta }(t,Y)=\left\vert P_{n}^{\prime }(t,Y)\right\vert \times \frac{%
1}{2\delta }1_{\{\left\vert y\right\vert \leq \delta \}}(P_{n}(t,Y)).
\label{d10}
\end{equation}%
Then, essentially, the K-R lemma says that for sufficiently small $\delta
_{n}$ we have 
\begin{eqnarray*}
\E(N_{n}(Y)) &\sim &\E(\int_{0}^{n\pi }\phi _{\delta _{n}}(t,Y)dt)\quad and \\
\E(N_{n}^{2}(Y)) &\sim &2\E(\int_{0}^{n\pi }dt\int_{0}^{t}\phi _{\delta
_{n}}(t,Y)\phi _{\delta _{n}}(s,Y)ds).
\end{eqnarray*}%
We make this precise in Lemma \ref{KR} bellow. Note that we will use the
above representations in connection with the $CLT$ - in particular we will
use the $CLT$ for $(\phi _{\delta _{n}}(t,Y),\phi _{\delta _{n}}(s,Y))$\ in
order to estimate $\E(\phi _{\delta _{n}}(t,Y)\phi _{\delta _{n}}(s,Y)).$ But
we will have to handle the following difficulty: if $t=s$ then the random
vector $(\phi _{\delta _{n}}(t,Y),\phi _{\delta _{n}}(s,Y))$ is degenerated,
so, in order to avoid this difficulty, we have to cancel a band around the
diagonal. The main ingredient in order to do it is the following lemma:

\begin{lemma}
\label{root} Let $I=(a,a+\varepsilon )$ and let $N_{n}(I,Y)=N_{a,a+%
\varepsilon }(P_{n}(.,Y))$ be the number of zeros of $P_{n}(t,Y)$ in $I.$
There exists universal constants $C\geq 1\geq c>0$ (independent of $%
n,a,\varepsilon )$ such that%
\begin{equation}
\E(N_{n}^{2}(I,Y)1_{\{N_{n}(I,Y)\geq 2\}})\leq C(\varepsilon ^{4/3}+ne^{-cn}).
\label{d2}
\end{equation}
\end{lemma}

\textbf{Proof.} Since the polynomial $P_{n}(t,Y)$ has at most $2n$ roots we
have%
\begin{equation}
\E(N_{n}^{2}(I,Y)1_{\{N_{n}(I,Y)\geq 2\}})=\P(N_{n}(Y)\geq
2)+\sum_{p=1}^{2n}(2p+1)\P(N_{n}(I,Y)>p)  \label{d3}
\end{equation}%
so we have to upper bound $\P(N_{n}(I,Y)>p).$ In order to do it we will use
the following fact: if $f:[a,a+\varepsilon ]\rightarrow R$ is $p+1$ times
differentiable and has at list $p+1$ zeros in this interval, then 
\begin{equation*}
\sup_{x\in \lbrack a,a+\varepsilon ]}\left\vert f(x)\right\vert \leq \frac{%
\varepsilon ^{p+1}}{(p+1)!}\sup_{x\in \lbrack a,a+\varepsilon ]}\left\vert
f^{(p+1)}(x)\right\vert .
\end{equation*}%
An argument which proves this is the following: Lagrange's interpolation
theorem says that given any $p+1$ points $x_{i},i=1,...,p+1$ in $%
[a,a+\varepsilon ]$ one may find a polynomial $P$ of order $p$ such that $%
P(x_{i})=f(x_{i})$ and $\sup_{x\in \lbrack a,a+\varepsilon ]}\left\vert
f(x)-P(x)\right\vert $ is upper bounded as in the previous inequality. Then
we take $x_{i},i=1,...,p+1$ to be the zeros of $f$ and, since $P$ is of
order $p$ and has $p+1$ roots, we have $P=0$ and we are done.

We denote $M_{n,p}=\sup_{t\in \lbrack a,a+\varepsilon ]}\left\vert
P_{n}^{(p+1)}(t,Y)\right\vert $ and we use the above inequality for $%
f(t)=P_{n}(t,Y)$ in order to obtain%
\begin{equation*}
\P(N_{n}(I,Y)>p)\leq \P(\left\vert P_{n}(a,Y)\right\vert \leq M\times \frac{%
(2\varepsilon )^{p+1}}{(p+1)!})+\P(M_{n,p}\geq M)
\end{equation*}%
A reasoning based on Sobolev's inequality and on Burkholder's inequality
(see the proof of Lemma 3.3 in the section ``small balls''of \cite{[BCP]})
proves that 
\begin{equation*}
\P(M_{n,p}\geq M)\leq \frac{1}{M^{2}}\E(M_{n,p}^{2})\leq \frac{C}{M^{2}}
\end{equation*}%
with $C$ a constant which depends on $p$ and on $M_{3}(Y)$.

We denote now $\delta =M\times \frac{(2\varepsilon )^{p+1}}{(p+1)!}$ and we
estimate%
\begin{equation*}
\P(\left\vert P_{n}(a,Y)\right\vert \leq \delta )=\delta \E(F_{\delta
}^{\prime }(P_{n}(a,Y)))\quad with\quad F_{\delta }(x)=\int_{-\infty
}^{x}1_{\{\left\vert y\right\vert \leq \delta \}}\frac{dy}{2\delta }.
\end{equation*}%
We will use (\ref{CLT}): we have $\E(\left\vert P_{n}(a,Y)\right\vert ^{2})=1$
and $\left\Vert F_{\delta }\right\Vert _{\infty }\leq 1$ and $\left\Vert
F_{\delta }^{\prime }\right\Vert _{\infty }\leq \delta ^{-1}$ so we get%
\begin{equation*}
\left\vert \E(F_{\delta }^{\prime }(P_{n}(a,Y))-\E(F_{\delta }^{\prime
}(W))\right\vert \leq C(\frac{1}{\sqrt{n}}+\frac{1}{\delta }e^{-cn})
\end{equation*}%
with $W$ a standard normal random variable. Since $\left\vert \E(F_{\delta
}^{\prime }(W))\right\vert \leq \frac{1}{2\pi }$ we get%
\begin{equation*}
\left\vert \E(F_{\delta }^{\prime }(P_{n}(a,Y))\right\vert \leq C(1+\frac{1}{%
\delta }e^{-cn}).
\end{equation*}%
This gives 
\begin{equation*}
\P(\left\vert P_{n}(a,Y)\right\vert \leq \delta )\leq C\delta +Ce^{-cn}
\end{equation*}%
and coming back%
\begin{equation*}
\P(N_{n}(I,Y)>p)\leq CM\times \frac{(2\varepsilon )^{p+1}}{(p+1)!}+Ce^{-cn}+%
\frac{C}{M^{2}}.
\end{equation*}%
We optimize on $M$ in order to obtain (for $p\geq 1)$%
\begin{equation*}
\P(N_{n}(I,Y)>p)\leq C\frac{\varepsilon ^{4/3}}{(p+1)!^{2/3}}+Ce^{-cn}.
\end{equation*}%
We insert this in (\ref{d3}) and, since $\sum_{p=1}^{\infty
}p/(p+1)!^{2/3}<\infty ,$ we obtain (\ref{d2}). $\square $

We fix now $\varepsilon >0,$ we denote 
\begin{equation*}
I_{k}^{\varepsilon }=[k\varepsilon ,(k+1)\varepsilon )\quad and\quad
D_{n,\varepsilon }=\cup _{0\leq k\leq n\pi /\varepsilon }\cup
_{p=0,k-2}I_{k}^{\varepsilon }\times I_{p}^{\varepsilon }
\end{equation*}%
We also denote%
\begin{eqnarray}
V_{n}(Y) &=&\E(N_{n}^{2}(Y))-(\E(N_{n}(Y)))^{2}\quad and  \label{VAR} \\
v_{n}(t,s,Y) &=&\E(\phi _{\delta _{n}}(t,Y)\phi _{\delta _{n}}(s,Y))-\E(\phi
_{\delta _{n}}(t,Y))\E(\phi _{\delta _{n}}(s,Y))  \notag
\end{eqnarray}%
with $N_{n}(Y)$ defined in (\ref{var}) and $\phi _{\delta _{n}}(t,Y)$
defined in (\ref{d10}).

\begin{lemma}
\label{KR}\textbf{A}. Let $\delta _{n}=n^{-\theta }$ with $\theta =5.$ Then 
\begin{equation}
\E(N_{n}^{2}(Y))=\E(N_{n}(Y))+2\int_{D_{n,\varepsilon }}\E(\phi _{\delta
_{n}}(t,Y)\phi _{\delta _{n}}(s,Y))dsdt+R_{n,\varepsilon }  \label{d8}
\end{equation}%
with 
\begin{equation}
\overline{\lim }_{n}\frac{1}{n}\left\vert R_{n,\varepsilon }\right\vert \leq
C\varepsilon ^{1/3}.  \label{d8''}
\end{equation}

\textbf{B}. And 
\begin{equation}
(\E(N_{n}(Y)))^{2}=2\int_{D_{n,\varepsilon }}\E(\phi _{\delta
_{n}}(t,Y))\E(\phi _{\delta _{n}}(s,Y))dsdt+R_{n,\varepsilon }  \label{d8'}
\end{equation}%
with $R_{n,\varepsilon }$ which verifies (\ref{d8''}).

\textbf{C}. 
\begin{equation}
V_{n}(Y)=V_{n}(G)+2\int_{D_{n,\varepsilon
}}(v_{n}(t,s,Y)-v_{n}(t,s,G))dsdt+R_{n,\varepsilon }  \label{d9}
\end{equation}%
with $R_{n,\varepsilon }$ which verifies (\ref{d8''}).
\end{lemma}

\textbf{Proof of A.} \textbf{Step 1}. We write%
\begin{equation*}
\E(N_{n}^{2}(Y))=J_{1}(n)+2J_{2}(n)+2J_{3}(n)
\end{equation*}%
with%
\begin{eqnarray*}
J_{1}(n) &=&\sum_{0\leq k\leq n\pi /\varepsilon
}\E(N_{n}^{2}(I_{k}^{\varepsilon },Y)),\quad J_{2}(n)=\sum_{0\leq k\leq n\pi
/\varepsilon }\E(N_{n}(I_{k}^{\varepsilon },Y)N_{n}(I_{k+1}^{\varepsilon },Y))
\\
J_{3}(n) &=&\sum_{0\leq k\leq n\pi /\varepsilon }\sum_{p=k+2}^{[n\pi
/\varepsilon ]}\E(N_{n}(I_{k}^{\varepsilon },Y)N_{n}(I_{p}^{\varepsilon },Y)).
\end{eqnarray*}%
Note that%
\begin{equation*}
\E(N_{n}(I_{k}^{\varepsilon },Y)N_{n}(I_{k+1}^{\varepsilon },Y))\leq
\E(N_{n}^{2}(I_{k}^{\varepsilon }\cup I_{k+1}^{\varepsilon
},Y)1_{\{N_{n}(I_{k}^{\varepsilon }\cup I_{k+1}^{\varepsilon },Y)\geq 2\}})
\end{equation*}%
Using (\ref{d2})%
\begin{equation*}
\left\vert J_{2}(n)\right\vert \leq C\times \frac{n}{\varepsilon }\times
(\varepsilon ^{4/3}+ne^{-n})
\end{equation*}%
so we get 
\begin{equation*}
\overline{\lim_{n}}\frac{1}{n}\left\vert J_{2}(n)\right\vert \leq
C\varepsilon ^{1/3}.
\end{equation*}%
We also have 
\begin{equation*}
\E(N_{n}^{2}(I_{k}^{\varepsilon },Y))=\E((N_{n}^{2}(I_{k}^{\varepsilon
},Y)-N_{n}(I_{k}^{\varepsilon },Y))1_{\{N_{n}(I_{k}^{\varepsilon },Y)\geq
2\}}))+\E(N_{n}(I_{k}^{\varepsilon },Y))
\end{equation*}%
so using (\ref{d2}) again 
\begin{equation*}
\overline{\lim_{n}}\frac{1}{n}\left\vert J_{n}(1)-\E(N_{n}(Y))\right\vert
\leq C\varepsilon ^{1/3}.
\end{equation*}

\textbf{Step 2}. We want to estimate 
\begin{equation*}
\frac{1}{n}J_{n}(3)=\frac{1}{n}\E(\sum_{0\leq k\leq n\pi /\varepsilon
}\sum_{p=k+2}^{[n\pi /\varepsilon ]}N_{n}(I_{k}^{\varepsilon
},Y)N_{n}(I_{p}^{\varepsilon },Y)).
\end{equation*}%
We will use the Kac-Rice formula (see the beginning of this section) for $%
f(t)=P_{n}(t,Y)$ so we have $N_{n}(Y)=N_{0,n\pi }(P_{n}(t,Y)).$ We denote $%
\delta _{n}(Y)=\delta _{0,n\pi }(P_{n}(.,Y)))$ (see (\ref{d3c}))$,$ we take $%
\delta _{n}=n^{-\theta }=n^{-5}$ and we write%
\begin{equation*}
\E(N_{n}(I_{k}^{\varepsilon },Y)N_{n}(I_{p}^{\varepsilon
},Y))=A_{n,k,p,\varepsilon }+B_{n,k,p,\varepsilon }
\end{equation*}%
with%
\begin{eqnarray*}
A_{n,k,p,\varepsilon } &=&\E(N_{n}(I_{k}^{\varepsilon
},Y)N_{n}(I_{p}^{\varepsilon },Y)1_{\{\delta _{n}\leq \delta _{n}(Y)\}}) \\
B_{n,k,p,\varepsilon } &=&\E(N_{n}(I_{k}^{\varepsilon
},Y)N_{n}(I_{p}^{\varepsilon },Y))1_{\{\delta _{n}>\delta _{n}(Y)\}}).
\end{eqnarray*}%
Since $P_{n}(t,Y)$ has at most $2n$ roots we get 
\begin{equation*}
B_{n,k,p,\varepsilon }\leq 4n^{2}\P(\delta _{n}\geq \delta _{n}(Y)).
\end{equation*}%
Recall that $\delta _{n}(Y)=\min \{\left\vert P_{n}(0,Y)\right\vert
,\left\vert P_{n}(n\pi ,Y)\right\vert ,\omega _{0,\pi }(P_{n})\}$ with $%
\omega _{0,\pi }(P_{n})=\inf_{0\leq t\leq n\pi }(\left\vert
P_{n}(t,Y)\right\vert +\left\vert P_{n}^{\prime }(t,Y)\right\vert ).$ Since $%
\left\vert P_{n}(0,Y)\right\vert =\left\vert P_{n}(n\pi ,Y)\right\vert =%
\frac{1}{\sqrt{n}}\left\vert \sum_{k=1}^{n}Y_{k}^{1}\right\vert $ it follows
that 
\begin{equation*}
\P(\delta _{n}\geq \delta _{n}(Y))\leq \P(\delta _{n}\geq \frac{1}{\sqrt{n}}%
\left\vert \sum_{k=1}^{n}Y_{k}^{1}\right\vert )+\P(\delta _{n}\geq \omega
_{0,\pi }(P_{n}))\leq \frac{C}{n^{4-\varepsilon }}
\end{equation*}%
the last inequality being a consequence of (\ref{SBB}) and (\ref{SB}) with $%
\theta =5.$ So we get 
\begin{equation*}
\frac{1}{n}\sum_{0\leq k\leq n\pi /\varepsilon }\sum_{p=k+2}^{[n\pi
/\varepsilon ]}B_{n,k,p,\varepsilon }\leq \frac{C}{n^{1-\varepsilon }}%
\rightarrow 0.
\end{equation*}

Moreover using K-R lemma (notice that $\delta _{n}(Y)\leq \delta
_{k\varepsilon ,(k+1)\varepsilon }(P_{n}(.,Y))$ for every $k)$\ we have%
\begin{equation*}
A_{n,k,p,\varepsilon }=\E(1_{\{\delta _{n}\leq \delta
_{n}(Y)\}}\int_{I_{k}^{\varepsilon }\times I_{p}^{\varepsilon }}\phi
_{\delta _{n}}(t,Y)\phi _{\delta _{n}}(s,Y)dtds)
\end{equation*}%
and consequently 
\begin{equation*}
\sum_{k=1}^{n\pi }\sum_{p=k+2}^{n\pi }A_{n,k,p,\varepsilon }=\E(1_{\{\delta
_{n}\leq \delta _{n}(Y)\}}\int_{D_{n,\varepsilon }}\phi _{\delta
_{n}}(t,Y)\phi _{\delta _{n}}(s,Y)dtds)=a_{n,\varepsilon }+b_{n,\varepsilon }
\end{equation*}%
with%
\begin{eqnarray*}
a_{n,\varepsilon } &=&\E(\int_{D_{n,\varepsilon }}\phi _{\delta
_{n}}(t,Y)\phi _{\delta _{n}}(s,Y)dtds) \\
b_{n,\varepsilon } &=&\E(1_{\{\delta _{n}\geq \delta
_{n}(Y)\}}\int_{D_{n,\varepsilon }}\phi _{\delta _{n}}(t,Y)\phi _{\delta
_{n}}(s,Y)dtds).
\end{eqnarray*}%
By (\ref{d3b}) 
\begin{eqnarray*}
\int_{D_{n,\varepsilon }}\phi _{\delta _{n}}(t,Y)\phi _{\delta
_{n}}(s,Y)dtds &\leq &(\int_{0}^{n\pi }\phi _{\delta _{n}}(t,Y)dt)^{2} \\
&\leq &(1+N_{n}([0,n\pi ],P_{n}^{\prime }(.,Y)))^{2}.
\end{eqnarray*}%
Since $P_{n}^{\prime }$ is still a trigonometric polynomial of order $n,$ it
has at most $2n$ roots. Then the above quantity is upper bounded by $%
(1+2n)^{2}$ and finally, using the small balls result%
\begin{equation*}
\frac{1}{n}b_{n,\varepsilon }\leq Cn^{3}\P(\delta _{n}\geq \delta
_{n}(Y))\leq \frac{C}{n^{1-\varepsilon }}\rightarrow 0
\end{equation*}%
so (\ref{d8}) is proved.

\textbf{Proof of B.} The proof is analogous (but simpler) so we just sketch
it. We denote by $R_{n}$ a quantity such that $\overline{\lim }_{n}\frac{1}{n%
}\left\vert R_{n}\right\vert =0.$\textbf{\ }Using again K-R formula and the
small balls property 
\begin{eqnarray*}
(\E(N_{n}(Y)))^{2} &=&(\E(\int_{0}^{n\pi }1_{\{\delta _{n}\leq \delta
_{n}(Y)\}}\phi _{\delta _{n}}(t,Y)dt)^{2}+R_{n} \\
&=&2\int_{0}^{n\pi }dt\int_{0}^{t}\E(1_{\{\delta _{n}\leq \delta
_{n}(Y)\}}\phi _{\delta _{n}}(t,Y))\E(1_{\{\delta _{n}\leq \delta
_{n}(Y)\}}\phi _{\delta _{n}}(s,Y))ds+R_{n} \\
&=&2\int_{0}^{n\pi }dt\int_{0}^{t}\E(\phi _{\delta _{n}}(t,Y))\E(\phi _{\delta
_{n}}(s,Y))ds+R_{n}^{\prime } \\
&=&2\int_{D_{n,\varepsilon }}\E(\phi _{\delta _{n}}(t,Y))\E(\phi _{\delta
_{n}}(s,Y))dsdt+R_{n,\varepsilon }+R_{n}^{\prime }
\end{eqnarray*}%
with%
\begin{equation*}
R_{n,\varepsilon }=\int_{D_{n,\varepsilon }^{c}}\E(\phi _{\delta
_{n}}(t,Y))\E(\phi _{\delta _{n}}(s,Y))dsdt.
\end{equation*}%
Using the $CLT$ we get%
\begin{equation*}
\left\vert \E(\phi _{\delta _{n}}(t,Y))-\E(\phi _{\delta
_{n}}(t,G))\right\vert \leq C(\frac{1}{\sqrt{n}}+\delta _{n}^{-1}e^{-cn}).
\end{equation*}%
Recall that $(P_{n}(t,G),P_{n}^{\prime }(t,G))$ is a Gaussian random
variable of covariance matrix $C_{n}(k,t)$ and, for sufficiently large $n$
one has $\left\langle C_{n}(k,t)x,x\right\rangle \geq \frac{1}{3}\left\vert
x\right\vert ^{2}.$ It follows that%
\begin{eqnarray*}
\E(\phi _{\delta _{n}}(t,G)) &=&\int_{R^{2}}\left\vert x_{2}\right\vert \frac{%
1}{2\delta _{n}}1_{\{\left\vert x_{1}\right\vert \leq \delta _{n}\}}\frac{1}{%
2\pi }e^{-\left\langle C_{n}(k,t)x,x\right\rangle }dx \\
&\leq &\int_{R}\left\vert x_{2}\right\vert \frac{1}{\sqrt{2\pi }}e^{-\frac{1%
}{6}\left\vert x_{2}\right\vert ^{2}}dx_{2}\times \int_{R}\frac{1}{2\delta
_{n}}1_{\{\left\vert x_{1}\right\vert \leq \delta _{n}\}}\frac{1}{\sqrt{2\pi 
}}e^{-\frac{1}{6}\left\vert x_{1}\right\vert ^{2}}dx_{1} \\
&\leq &C.
\end{eqnarray*}%
So $\E(\phi _{\delta _{n}}(t,Y))\leq C$ and consequently, for sufficiently
large $n$%
\begin{equation*}
\frac{1}{n}\left\vert R_{n,\varepsilon }\right\vert \leq \frac{C}{n}%
\left\vert D_{n,\varepsilon }^{c}\right\vert \leq C\varepsilon .
\end{equation*}

\textbf{Proof of C.} We have proved in \cite{[BCP]} that 
\begin{equation*}
\lim_{n}\frac{1}{n}(\E(N_{n}(Y))-\E(N_{n}(G)))=0
\end{equation*}%
so (\ref{d9}) is an immediate consequence of (\ref{d8}) and (\ref{d8'}). $%
\square $

\section{Cancellations}

\label{Cancellations}

Having in mind (\ref{d9})\ we will now estimate $v_{n}(t,s,Y)-v_{n}(t,s,G).$
A careful analysis of this term involve a certain number of cancellations.
The objects which are involved here are the following. For each $t\geq 0$ we
consider the matrices $C_{n}(r,t)\in \mathcal{M}_{2\times 2},n\in N,1\leq
r\leq n$ defined by%
\begin{equation}
C_{n}(r,t)=\left( 
\begin{tabular}{ll}
$\cos (\frac{rt}{n})$ & $\sin (\frac{rt}{n})$ \\ 
$-\frac{r}{n}\sin (\frac{rt}{n})$ & $\frac{r}{n}\cos (\frac{rt}{n})$%
\end{tabular}%
\right) .  \label{N1}
\end{equation}%
Moreover we consider the sequence $Y=(Y_{r})_{r\in N}$ introduced in the
first section and we denote 
\begin{equation}
Z_{n,r}(t,Y)=C_{n}(r,t)Y_{r}.  \label{N2}
\end{equation}%
We are concerned with%
\begin{equation}
S_{n}(t,Y)=\frac{1}{\sqrt{n}}\sum_{r=1}^{n}Z_{n,r}(t,Y)=\frac{1}{\sqrt{n}}%
\sum_{r=1}^{n}C_{n}(r,t)Y_{r}.  \label{N3}
\end{equation}%
Notice that, with the notation form (\ref{pol}), \ $%
S_{n}^{1}(t,Y)=P_{n}(t,Y) $ and $S_{n}^{2}(t,Y)=P_{n}^{\prime }(t,Y).$ We
also denote%
\begin{equation}
S_{n}(t,s,Y)=\left( 
\begin{tabular}{l}
$S_{n}(t,Y)$ \\ 
$S_{n}(s,Y)$%
\end{tabular}%
\right) ,\quad Z_{n,r}(t,s,Y)=\left( 
\begin{tabular}{l}
$Z_{n,r}(t,Y)$ \\ 
$Z_{n,r}(s,Y)$%
\end{tabular}%
\right) .  \label{N4}
\end{equation}%
In order to be able to give our results we need to introduce some more
notation. We denote 
\begin{equation}
\Phi _{\delta }(x_{1},x_{2})=\left\vert x_{2}\right\vert f_{\delta
}(x_{1})\quad with\quad f_{\delta }(x_{1})=\frac{1}{2\delta }\int_{-\infty
}^{x_{1}}1_{\{\left\vert y\right\vert \leq \delta \}}dy  \label{N8}
\end{equation}%
so that $\partial _{1}\Phi _{\delta }(S_{n}(t,Y))=\phi _{\delta }(t,Y)$ (see
(\ref{d10})). Then (see (\ref{VAR}) and recall that $\delta _{n}=1/n^{5}$) 
\begin{equation}
v_{n}(t,s,Y)=\E(\partial _{1}\Phi _{\delta _{n}}(S_{n}(t,Y))\partial _{1}\Phi
_{\delta _{n}}(S_{n}(s,Y)))-\E(\partial _{1}\Phi _{\delta
_{n}}(S_{n}(t,Y)))\E(\partial _{1}\Phi _{\delta _{n}}(S_{n}(s,Y)))
\label{N8'}
\end{equation}%

\begin{lemma}
\label{Cancel}Suppose that\ for every multi-index $\alpha $ with $\left\vert
\alpha \right\vert =3,4$\ the following limits exists and are finite:%
\begin{equation*}
\lim_{n}\E(\prod_{i=1}^{\left\vert \alpha \right\vert }Y_{n}^{\alpha
_{i}})=y_{\infty }(\alpha ).
\end{equation*}%
Then, for every $\varepsilon >0,$ 
\begin{equation}
\lim_{n}\frac{1}{n}\int_{D_{n,\varepsilon
}}(v_{n}(t,s,Y)-v_{n}(t,s,G))dsdt= \frac{1}{120}\times  y_{\ast }+r_{\varepsilon }
\label{N4'}
\end{equation}%
with $\left\vert r_{\varepsilon }\right\vert \leq C\varepsilon $ and 
\begin{equation*}
y_{\ast }=(y_{\infty }(1,1,2,2)-1)+(y_{\infty }(2,2,1,1)-1)+(y_{\infty
}(1,1,1,1)-3)+(y_{\infty }(2,2,2,2)-3).
\end{equation*}
\end{lemma}

\bigskip

\textbf{Proof. Step 1.} We come back to the framework from Section 2. We
denote by $\Sigma _{n}(t)$ the covariance matrix of $S_{n}(t,Y)$ and by $%
\Sigma _{n}(t,s)$ the covariance matrix of $S_{n}(t,s,Y)$ and we will use (%
\ref{N7e}) for $X_{n,k}=Z_{n,k}(t,Y)$ respectively for $%
X_{n,k}=Z_{n,k}(t,s,Y)=(Z_{n,k}(t,Y),Z_{n,k}(s,Y)).$ We stress that all the
constants will depend on $\det \Sigma _{n}(t,s)$ which is larger then $\frac{%
1}{2}\lambda ^{2}(\varepsilon )>0$ for $(t,s)\in D_{n,\varepsilon }$ (see (%
\ref{App5}))$.$ We will also use the diagonal matrices $I_{2}=I_{2}(\lambda
) $ with $\lambda _{1}=1,\lambda _{2}=\frac{1}{3}$ and $I_{4}=I_{4}(\lambda
) $ with $\lambda _{1}=\lambda _{3}=1,\lambda _{2}=\lambda _{4}=\frac{1}{3}.$
By (\ref{N7e})%
\begin{eqnarray}
\E(\partial _{1}\Phi _{\delta _{n}}(S_{n}(t,Y))) &=&\E(\partial _{1}\Phi
_{\delta _{n}}(S_{n}(t,G)))  \label{R1} \\
&&+\frac{1}{n}\E(\partial _{1}\Phi _{\delta _{n}}(I_{2}^{1/2}W)\Gamma
_{n,2}(I_{2}^{-1/2}Z_{n}(t,Y),W)))  \notag \\
&&+\frac{1}{n}r_{n}(t,\Phi _{\delta _{n}})+\frac{1}{n^{3/2}}R_{n}(t,\Phi
_{\delta _{n}})  \notag
\end{eqnarray}%
and a similar expression for $\E(\partial _{1}\Phi _{\delta
_{n}}(S_{n}(s,Y))).$ The remainder $r_{n}(t,\Phi _{\delta })$ verifies (\ref%
{N7d}) with $\Sigma _{n}(t)-I_{2}(\lambda ).$ We also recall that $%
S_{n}(t,G) $ has the same law as $\Sigma _{n}^{1/2}(t)W$ so, (with $%
r_{n}(t,\Phi _{\delta })$ which verifies (\ref{N7d})), 
\begin{equation*}
\E(\partial _{1}\Phi _{\delta _{n}}(S_{n}(t,G)))=\E(\partial _{1}\Phi _{\delta
_{n}}(\Sigma _{n}^{1/2}(t)W))=\E(\partial _{1}\Phi _{\delta
_{n}}(I_{2}^{1/2}W))+\frac{1}{n}r_{n}(t,\Phi _{\delta }).
\end{equation*}

Moreover, we denote $\Psi _{\delta }(x_{1},x_{2},x_{3},x_{4})=\Phi _{\delta
}(x_{1},x_{2})\Phi _{\delta }(x_{3},x_{4})$ and we write%
\begin{eqnarray}
\E(\partial _{1}\partial _{3}\Psi _{\delta _{n}}(S_{n}(t,s,Y))) &=&\E(\partial
_{1}\partial _{3}\Psi _{\delta _{n}}(S_{n}(t,s,G)))  \label{R2} \\
&&+\frac{1}{n}\E(\partial _{1}\partial _{3}\Psi _{\delta
_{n}}(I_{4}^{1/2}W)\Gamma _{n,2}(I_{4}^{-1/2}Z_{n}(t,s,Y),W)))  \notag \\
&&+\frac{1}{n}r_{n}(t,s,\Psi _{\delta _{n}})+\frac{1}{n^{3/2}}R_{n}(t,s,\Psi
_{\delta _{n}}).  \notag
\end{eqnarray}%
Here \ $r_{n}(t,s,\Psi _{\delta _{n}})$ verifies (\ref{N7d}) with $\Sigma
_{n}(t,s)-I_{4}(\lambda ).$ And, as above,%
\begin{equation*}
\E(\partial _{1}\partial _{3}\Psi _{\delta _{n}}(S_{n}(t,s,G)))=\E(\partial
_{1}\partial _{3}\Psi _{\delta _{n}}(\Sigma _{n}^{1/2}(t,s)W))=\E(\partial
_{1}\partial _{3}\Psi _{\delta _{n}}(I_{4}^{1/2}W))+\frac{1}{n}%
r_{n}(t,s,\Psi _{\delta }).
\end{equation*}

Our aim now is to estimate $v_{n}(t,s,Y)-v_{n}(t,s,G)$ (recall that $%
v_{n}(t,s,Y)$ is defined in (\ref{N8'}))$.$ In order to simplify notation we
put%
\begin{eqnarray*}
A_{n}(t,Y) &=&\E(\partial _{1}\Phi _{\delta _{n}}(S_{n}(t,Y))),\quad
A_{n}(t,s,Y)=\E(\partial _{1}\partial _{3}\Psi _{\delta _{n}}(S_{n}(t,s,Y))),
\\
C_{n}(t) &=&\E(\partial _{1}\Phi _{\delta _{n}}(I_{2}^{1/2}W)\Gamma
_{n,2}(I_{2}^{-1/2}Z_{n}(t,Y),W))), \\
C_{n}(t,s) &=&\E(\partial _{1}\partial _{3}\Psi _{\delta
_{n}}(I_{4}^{1/2}W)\Gamma _{n,2}(I_{4}^{-1/2}Z_{n}(t,s,Y),W))), \\
\widehat{R}_{n}(t) &=&\frac{1}{n}r_{n}(t,\Phi _{\delta _{n}})+\frac{1}{%
n^{3/2}}R_{n}(t,\Phi _{\delta _{n}}),\quad \widehat{R}_{n}(t,s)=\frac{1}{n}%
r_{n}(t,s,\Psi _{\delta _{n}})+\frac{1}{n^{3/2}}R_{n}(t,s,\Psi _{\delta
_{n}}).
\end{eqnarray*}%
With this notation (\ref{R1}) and (\ref{R2}) read%
\begin{eqnarray*}
A_{n}(t,Y) &=&A_{n}(t,G)+\frac{1}{n}C_{n}(t)+\widehat{R}_{n}(t), \\
A_{n}(t,s,Y) &=&A_{n}(t,s,G)+\frac{1}{n}C_{n}(t,s)+\widehat{R}_{n}(t,s)
\end{eqnarray*}%
and consequently%
\begin{equation*}
v_{n}(t,s,Y)-v_{n}(t,s,G)=\frac{1}{n}\gamma _{n}(t,s)+\overline{R}_{n}(t,s)
\end{equation*}%
with%
\begin{equation*}
\overline{R}_{n}(t,s)=(\frac{1}{n}C_{n}(t)+\widehat{R}_{n}(t))(\frac{1}{n}%
C_{n}(s)+\widehat{R}_{n}(s))-\widehat{R}_{n}(t)A_{n}(s,Y)-\widehat{R}%
_{n}(s)A_{n}(t,Y)
\end{equation*}%
and%
\begin{eqnarray*}
\gamma _{n}(t,s) &=&C_{n}(t,s)-C_{n}(t)A_{n}(s,Y)-C_{n}(s)A_{n}(t,Y) \\
&=&\E(\partial _{1}\partial _{3}\Psi _{\delta _{n}}(I_{4}^{1/2}W)\Gamma
_{n,2}(I_{4}^{-1/2}Z_{n}(t,s,Y),W))) \\
&&-\E(\partial _{1}\Phi _{\delta _{n}}(I_{2}^{1/2}W))\times \E(\partial
_{1}\Phi _{\delta _{n}}(I_{2}^{1/2}W)\Gamma
_{n,2}(I_{2}^{-1/2}Z_{n}(s,Y),W))) \\
&&-\E(\partial _{1}\Phi _{\delta _{n}}(I_{2}^{1/2}W)\Gamma
_{n,2}(I_{2}^{-1/2}Z_{n}(t,Y),W)))\times \E(\partial _{1}\Phi _{\delta
_{n}}(I_{2}^{1/2}W)).
\end{eqnarray*}

Notice that in the above expression of $\gamma _{n}(t,s)$, $W$ stands for a
standard normal random variable which is in dimension $4$ in the first
expectation and in dimension two in the following two ounces. In order to
put everything together we take two independent two-dimensional standard
normal random variables $W^{\prime }$ and $W^{\prime \prime }$ and we put $%
W=(W^{\prime },W^{\prime \prime })\in R^{4}$ which is itself a standard
normal random variable. Then 
\begin{equation*}
\partial _{1}\Phi _{\delta _{n}}(I_{2}^{1/2}W^{\prime })\partial _{1}\Phi
_{\delta _{n}}(I_{2}^{1/2}W^{\prime \prime })=\partial _{1}\partial _{3}\Psi
_{\delta _{n}}(I_{4}^{1/2}W)
\end{equation*}%
so we obtain 
\begin{eqnarray*}
\gamma _{n}(t,s) &=&\E(\partial _{1}\partial _{3}\Psi _{\delta
_{n}}(I_{4}^{1/2}W)[\Gamma _{n,2}(I_{4}^{-1/2}Z_{n}(t,s,Y),W))) \\
&&-\Gamma _{n,2}(I_{2}^{-1/2}Z_{n}(t,Y),W^{\prime })-\Gamma
_{n,2}(I_{2}^{-1/2}Z_{n}(s,Y),W^{\prime \prime })]).
\end{eqnarray*}%
We recall the definitions of $\Gamma _{n,2}^{\prime },\Gamma _{n,2}^{\prime
\prime }$ given in (\ref{N7a}) and we write $\gamma _{n}(t,s)=\gamma
_{n}^{\prime }(t,s)+\gamma _{n}^{\prime \prime }(t,s)$ with $\gamma ^{\prime
}$ which involves $\Gamma ^{\prime }$ and $\gamma ^{\prime \prime }$ which
involves $\Gamma ^{\prime \prime }$ instead of $\Gamma .$ We will analyze
them separately.

\textbf{Step 2}. Estimate of $\gamma ^{\prime \prime }.$ Our aim is to prove
that%
\begin{equation}
\frac{1}{n^{2}}\int_{D_{n,\varepsilon }}\gamma _{n}^{\prime \prime
}(t,s)dsdt=\int_{0}^{\pi }\int_{0}^{\pi }1_{D_{n,\varepsilon }}(nt,ns)\gamma
_{n}^{\prime \prime }(nt,ns)dsdt\rightarrow 0.  \label{N8''}
\end{equation}%
The analysis is based on (\ref{N7aa}). There are two kinds of cancellation
which are at work:

\textbf{First cancellation (mixed multi-indexes).} Denote $m_{k}(I)$\ the
set of the multi-indexes $\alpha =(\alpha _{1},...,\alpha _{k})$ with $%
\alpha _{i}\in I.$ Recall that $W=(W^{\prime },W^{\prime \prime })$ and
notice that if $\alpha \in m_{3}(1,2)$ then $H_{\alpha }(W)=H_{\alpha
}(W^{\prime }).$ But, if $\alpha \in m_{3}(3,4),$ then one has $H_{\alpha
}(W)=H_{\alpha }(W^{3},W^{4})=H_{\alpha }((W^{\prime \prime
})^{1},(W^{\prime \prime })^{2}).$ This means that, in the second case, a
``change of variable'' is needed: $\alpha =(\alpha _{1},\alpha _{2},\alpha
_{3})\rightarrow \widehat{\alpha }=(\alpha _{1}-2,\alpha _{2}-2,\alpha
_{3}-2)$: for example $(3,3,4)\rightarrow (1,1,2)$ or $(4,4,3)\rightarrow
(2,2,1).$ Having this in mind we go on and analyze $\Gamma _{n,2}^{\prime
\prime }$ defined in (\ref{N7aa}):\ 

\begin{eqnarray*}
\Gamma _{n,2}^{\prime \prime }(I_{4}^{-1/2}Z_{n}(t,s,Y),W)) &=&\frac{1}{72}%
\sum_{\left\vert \rho \right\vert =3}\sum_{\left\vert \beta \right\vert
=3}c_{n}(\beta ,I_{4}^{-1/2}Z_{n}(t,s,Y))c_{n}(\rho
,I_{4}^{-1/2}Z_{n}(t,s,Y))H_{(\beta ,\rho )}(W) \\
\Gamma _{n,2}^{\prime \prime }(I_{4}^{-1/2}Z_{n}(t,Y),W)) &=&\frac{1}{72}%
\sum_{\left\vert \rho \right\vert =3}\sum_{\left\vert \beta \right\vert
=3}c_{n}(\beta ,I_{4}^{-1/2}Z_{n}(t,Y))c_{n}(\rho
,I_{4}^{-1/2}Z_{n}(t,Y))H_{(\beta ,\rho )}(W^{\prime }) \\
\Gamma _{n,2}^{\prime \prime }(I_{4}^{-1/2}Z_{n}(s,Y),W)) &=&\frac{1}{72}%
\sum_{\left\vert \rho \right\vert =3}\sum_{\left\vert \beta \right\vert
=3}c_{n}(\beta ,I_{4}^{-1/2}Z_{n}(s,Y))c_{n}(\rho
,I_{4}^{-1/2}Z_{n}(s,Y))H_{(\beta ,\rho )}(W^{\prime \prime }).
\end{eqnarray*}%
Notice that the multi-indexes in the first line belong to $m_{3}(1,2,3,4)$
while the multi-indexes in the second and in the third line belong to $%
m_{3}(1,2)$. We look now to the sums in the first line. If all the elements
of $(\beta ,\rho )$ belong to $\{1,2\}$ then $H_{(\beta ,\rho
)}(W)=H_{(\beta ,\rho )}(W^{\prime })$ and $c_{n}(\beta
,I_{4}^{-1/2}Z_{n}(nt,ns,Y))=c_{n}(\beta ,I_{2}^{-1/2}Z_{n}(nt,Y))$ so the
corresponding term cancels. In the same way, if all the elements of $(\beta
,\rho )$ belong to $\{3,4\}$ then $H_{(\beta ,\rho )}(W)=H_{(\widehat{\beta }%
,\widehat{\rho })}(W^{\prime \prime })$ and $c_{n}(\beta
,I_{4}^{-1/2}Z_{n}(nt,ns,Y))=c_{n}(\widehat{\beta },I_{2}^{-1/2}Z_{n}(ns,Y))$
and the corresponding term cancels as well. We remain with ``mixed
multi-indexes'', such that $(\beta ,\rho )$ contain at least one element
from each of $\{1,2\}$ and of $\{3,4\}$.

\textbf{Second cancellation (even multi-indexes). }For each $i=1,...,4$ the
function $W_{i}\rightarrow \partial _{1}\partial _{3}\Psi _{\delta
_{n}}(I_{4}^{1/2}W)$ is even, so, because the symmetry argument 
\begin{equation*}
\E(\partial _{1}\partial _{3}\Psi _{\delta _{n}}(I_{4}^{1/2}W))H_{(\rho
,\beta )}(W))=0
\end{equation*}%
except the case when all the elements in $(\rho ,\beta )$ appear an even
number of times (this means that $i_{j}((\rho ,\beta ))$ is even for every $%
j=1,...,4).$

There are three types of multi-indexes which verify both conditions: take $%
i\in \{1,2\}$ and $j,p\in \{3,4\}$ (or the converse).%
\begin{eqnarray}
Case 1\quad \rho &=&(i,j,j),\quad \beta =(i,p,p)  \label{N8a} \\
Case 2\quad \rho &=&(i,i,j),\quad \beta =(j,p,p)  \label{N8b} \\
Case 3\quad \rho &=&(i,j,p),\quad \beta =(i,j,p)  \label{N8c}
\end{eqnarray}

We treat the Case 1 (the other cases are similar). In order to fix the ideas
we take $i=1$ and $j=4,$ so that $\rho =(1,4,4)$ (all the other cases are
similar). We compute%
\begin{equation*}
\E((C_{n}(k,nt)Y_{k})^{1}((C_{n}(k,ns)Y_{k})^{2})^{2})=%
\sum_{l_{1},l_{2},l_{3}=1}^{2}C_{n}^{1,l_{1}}(k,nt)C_{n}^{2,l_{2}}(k,ns)C_{n}^{2,l_{3}}(k,ns)\E(\prod_{i=1}^{3}Y_{k}^{l_{i}})
\end{equation*}%
Since in the Gaussian case we have $\E(\prod_{i=1}^{3}G_{k}^{l_{i}})=0,$ we
conclude that 
\begin{equation*}
\Delta _{\rho
}(Z_{n,k}(nt,ns,Y))=%
\sum_{l_{1},l_{2},l_{3}=1}^{2}C_{n}^{1,l_{1}}(k,nt)C_{n}^{2,l_{2}}(k,s)C_{n}^{2,l_{3}}(k,ns)\E(\prod_{i=1}^{3}Y_{k}^{l_{i}})
\end{equation*}%
and then 
\begin{equation*}
c_{n}(\rho ,I_{4}^{-1/2}Z_{n}(nt,ns,Y))=c_{n}^{\prime }(\rho
,I_{4}^{-1/2}Z_{n}(nt,ns,Y))+c_{n}^{\prime \prime }(\rho
,I_{4}^{-1/2}Z_{n}(nt,ns,Y))
\end{equation*}%
with%
\begin{eqnarray*}
c_{n}^{\prime }(\rho ,I_{4}^{-1/2}Z_{n}(nt,ns,Y))
&=&\sum_{l_{1},l_{2},l_{3}=1}^{2}y_{\infty }(l_{1},l_{2},l_{3})\times \frac{1%
}{n}%
\sum_{k=1}^{n}C_{n}^{1,l_{1}}(k,nt)C_{n}^{2,l_{2}}(k,s)C_{n}^{2,l_{3}}(k,ns)
\\
c_{n}^{\prime \prime }(\rho ,I_{4}^{-1/2}Z_{n}(nt,ns,Y))
&=&\sum_{l_{1},l_{2},l_{3}=1}^{2}\frac{1}{n}%
\sum_{k=1}^{n}C_{n}^{1,l_{1}}(k,nt)C_{n}^{2,l_{2}}(k,s)C_{n}^{2,l_{3}}(k,ns)(\E(\prod_{i=1}^{3}Y_{k}^{l_{i}})-y_{\infty }(l_{1},l_{2},l_{3})).
\end{eqnarray*}%
Since $\left\vert C_{n}^{i,j}(k,u)\right\vert \leq 1$ for every $i,j\in
\{1,2\}$ and $u>0,$ we have 
\begin{equation*}
\left\vert c_{n}^{\prime \prime }(\rho
,I_{4}^{-1/2}Z_{n}(nt,ns,Y))\right\vert \leq \sum_{l_{1},l_{2},l_{3}=1}^{2}%
\frac{1}{n}\sum_{k=1}^{n}\left\vert
\E(\prod_{i=1}^{3}Y_{k}^{l_{i}})-y_{\infty }(l_{1},l_{2},l_{3})\right\vert
\rightarrow 0.
\end{equation*}%
And using (\ref{er11}) we get $c_{n}^{\prime }(\rho
,I_{4}^{-1/2}Z_{n}(nt,ns,Y))\rightarrow 0.$ This is true for $t$ and $s$
such that $\frac{t}{\pi },\frac{s}{\pi },\frac{t+s}{\pi }$ and $\frac{t-s}{%
\pi }$ are irrational. But this means that this is true $dtds$ almost
surely. Then, using Lebesgue's dominated convergence theorem (notice that
the coefficients $c_{n},n\in N$ are uniformly bounded) we get 
\begin{equation*}
\int_{0}^{\pi }\int_{0}^{\pi }1_{D_{n,\varepsilon }}(nt,ns)c_{n}(\rho
,I_{4}^{-1/2}Z_{n}(nt,ns,Y))c_{n}(\beta
,I_{4}^{-1/2}Z_{n}(nt,ns,Y))dtds\rightarrow 0.
\end{equation*}%
So we have finished to prove (\ref{N8''}).

\textbf{Step 3}. We compute now%
{}{
\begin{equation*}
\lim_{n}\frac{1}{n^{2}}\int_{D_{n,\varepsilon }}\gamma _{n}^{\prime
}(t,s)dsdt=\frac{1}{2}\lim_{n}\frac{1}{n^{2}}\int_{0}^{n\pi }\int_{0}^{n\pi }\gamma
_{n}^{\prime }(t,s)dsdt+O(\epsilon)=\frac{1}{2}\lim_{n}\int_{0}^{\pi }\int_{0}^{\pi }\gamma
_{n}^{\prime }(nt,ns)dsdt+O(\epsilon),
\end{equation*}
where $O(\epsilon)$ is uniform in $n$.} We recall (\ref{N7b}). As in the previous discussion we notice that we have
two kind of cancellations: if all the components of $\alpha $ belong to $%
\{1,2\}$ or to $\{3,4\}$ then the corresponding term cancels. And for
symmetry reasons one also needs to have each component of $\alpha $ an even
number of times. So the only multi-indexes which have a non null
contribution are ({}{up to permutations}) $\alpha =(i,i,j,j)$ with $i\in \{1,2\}$ and $j\in \{3,4\}.$ {}{More precisely, for every fixed $(i,j)\in \{1,2\}\times\{3,4\}$ the following multi-indexes bring a non zero contribution: $(i,i,j,j), (i,j,i,j), (i,j,j,i), (j,j,i,i), (j,i,j,i),(j,i,i,j)$. Besides all the forthcoming computations are independent of the chosen permutations and we will simply assume that the multi-index is $(i,i,j,j)$ and multiply the final result by a factor $6$}. {}{Indeed, we} observe that%
\begin{equation*}
\gamma _{n}^{\prime }(nt,ns)=\frac{1}{24}\sum_{\alpha }\E(\partial
_{1}\partial _{3}\Psi _{\delta _{n}}(I_{4}^{1/2}W)H_{\alpha
}(W))c_{n}(\alpha ,I_{4}^{-1/2}Z_{n}(nt,ns,Y))
\end{equation*}%
with the sum over the multi-indexes of the form {}{(up to permutations)} $\alpha =(i,i,j,j)$ with $%
i\in \{1,2\}$ and $j\in \{3,4\}.$

We fix such a multi index $\alpha =(i,i,j,j)$ and we denote (with $j^{\prime
}=j-2)$%
\begin{eqnarray*}
p(\alpha ) &=&3^{i+j-4}=3^{i+j^{\prime}-2},\quad \\
U(\alpha ) &=&\frac{p(\alpha )y_{\ast }}{4(1+2(i+j-4))}=\frac{p(\alpha
)y_{\ast }}{4(1+2(i+j^{\prime }-2))}.
\end{eqnarray*}%
Our first aim is to prove that, if $\frac{t}{\pi },\frac{s}{\pi },\frac{t+s}{%
\pi },\frac{t-s}{\pi }$ are irrational, then 
\begin{equation}
\lim_{n}c_{n}(\alpha ,I_{4}^{-1/2}Z_{n}(nt,ns,Y)=U(\alpha ).  \label{lim}
\end{equation}%
We compute 
\begin{eqnarray*}
&&\E((C_{n}(k,nt)Y_{k})^{i})^{2}(C_{n}(k,ns)Y_{k})^{j-2})^{2}) \\
&=&%
\sum_{l_{1},l_{2},l_{3},l_{4}=1}^{2}C_{n}^{i,l_{1}}(k,nt)C_{n}^{i,l_{2}}(k,nt)C_{n}^{j-2,l_{3}}(k,ns)C_{n}^{j-2,l_{4}}(k,ns)\E(\prod_{i=1}^{4}Y_{k}^{l_{i}}).
\end{eqnarray*}%
Then%
\begin{eqnarray*}
&&\Delta _{\alpha }(I_{4}^{-1/2}Z_{n,k}(nt,ns,Y)) \\
&=&p(\alpha
)\sum_{l_{1},l_{2},l_{3},l_{4}=1}^{2}C_{n}^{i,l_{1}}(k,nt)C_{n}^{i,l_{2}}(k,nt)C_{n}^{j-2,l_{3}}(k,ns)C_{n}^{j-2,l_{4}}(k,ns)(\E(\prod_{i=1}^{4}Y_{k}^{l_{i}})-\E(\prod_{i=1}^{4}G_{k}^{l_{i}}))
\end{eqnarray*}%
and finally%
\begin{eqnarray*}
c_{n}(\alpha ,I_{4}^{-1/2}Z_{n}(nt,ns,Y)) &=&\frac{1}{n}\sum_{k=1}^{n}\Delta
_{\alpha }(I_{4}^{-1/2}Z_{n,k}(nt,ns,Y)) \\
&=&p(\alpha )\sum_{l_{1},l_{2},l_{3},l_{4}=1}^{2}c_{n}^{\prime }(\alpha
,l_{1},l_{2},l_{3},l_{4})+c_{n}^{\prime \prime }(\alpha
,l_{1},l_{2},l_{3},l_{4})
\end{eqnarray*}%
with%
\begin{eqnarray*}
c_{n}^{\prime }(\alpha ,l_{1},l_{2},l_{3},l_{4}) &=&\frac{1}{n}%
\sum_{k=1}^{n}C_{n}^{i,l_{1}}(k,nt)C_{n}^{i,l_{2}}(k,nt)C_{n}^{j-2,l_{3}}(k,ns)C_{n}^{j-2,l_{4}}(k,ns)(y_{\infty }(l_{1},l_{2},l_{3},l_{4})-\E(\prod_{i=1}^{4}B^{l_{i}}))
\\
c_{n}^{\prime \prime }(\alpha ,l_{1},l_{2},l_{3},l_{4}) &=&\frac{1}{n}%
\sum_{k=1}^{n}C_{n}^{i,l_{1}}(k,nt)C_{n}^{i,l_{2}}(k,nt)C_{n}^{j-2,l_{3}}(k,ns)C_{n}^{j-2,l_{4}}(k,ns)(\E(\prod_{i=1}^{4}Y_{k}^{l_{i}})-y_{\infty }(l_{1},l_{2},l_{3},l_{4})).
\end{eqnarray*}%
Here $B=(B^{1},B^{2})$ is a standard Gaussian random variable. Since $%
\E(\prod_{i=1}^{4}Y_{k}^{l_{i}})\rightarrow y_{\infty
}(l_{1},l_{2},l_{3},l_{4})$ we get $c_{n}^{\prime \prime }(\alpha
,l_{1},l_{2},l_{3},l_{4})\rightarrow 0.$ We analyze now $c_{n}^{\prime
}(\alpha ,l_{1},l_{2},l_{3},l_{4}).$ By (\ref{er10}), if $l_{1}\neq l_{2}$
or if $l_{3}\neq l_{4}$ this term converges to zero$.$ So we have to
consider only 
\begin{equation*}
c_{n}^{\prime }(\alpha ,l,l,l^{\prime },l^{\prime })=\frac{1}{n}%
\sum_{k=1}^{n}C_{n}^{i,l}(k,nt)^{2}C_{n}^{j-2,l^{\prime
}}(k,ns)^{2}(y_{\infty }(l,l,l^{\prime },l^{\prime
})-\E((B^{l})^{2}(B^{l^{\prime }})^{2})
\end{equation*}%
Take first $l=1$ and $l^{\prime }=2.$ Then, using (\ref{er9}), we have 
\begin{eqnarray*}
c_{n}^{\prime }(\alpha ,1,1,2,2) &=&\frac{1}{n}%
\sum_{k=1}^{n}C_{n}^{i,l}(k,nt)^{2}C_{n}^{j-2,l^{\prime
}}(k,ns)^{2}(y_{\infty }(1,1,2,2)-1) \\
&\rightarrow &\frac{1}{4(1+2(i+j-4))}(y_{\infty }(1,1,2,2)-1).
\end{eqnarray*}%
And if $l=l^{\prime }=1$ (or if $l=l^{\prime }=2)$ we have 
\begin{eqnarray*}
c_{n}^{\prime }(\alpha ,1,1,1,1) &=&\frac{1}{n}%
\sum_{k=1}^{n}C_{n}^{i,l}(k,nt)^{2}C_{n}^{j-2,l^{\prime
}}(k,ns)^{2}(y_{\infty }(1,1,1,1)-3) \\
&\rightarrow &\frac{1}{4(1+2(i+j-4))}(y_{\infty }(1,1,1,1)-3).
\end{eqnarray*}%
So (\ref{lim}) is proved and, as an immediate consequence we obtain%
\begin{equation}
\lim_{n}\int_{0}^{\pi }\int_{0}^{\pi }c_{n}(\alpha
,I_{4}^{-1/2}Z_{n}(nt,ns,Y))dsdt=\pi ^{2}U(\alpha ).  \label{lim1}
\end{equation}

We compute now%
\begin{equation*}
\lim_{n}\E(\partial _{1}\partial _{3}\Psi _{\delta
_{n}}(I_{4}^{1/2}W)H_{\alpha }(W)).
\end{equation*}%
Notice that if $i\in \{1,2\}$ and $j\in \{3,4\}$ then (recall that $h_{2}$
is the Hermite polynomial of order $2$ on $R,$ so that $h_{2}(x)=x^{2}-1)$%
\begin{equation*}
H_{(i,i,j,j)}(W)=h_{2}(W_{i}^{\prime })h_{2}(W_{j-2}^{\prime \prime }))
\end{equation*}%
so that%
\begin{eqnarray*}
\E(\partial _{1}\partial _{3}\Psi _{\delta _{n}}(I_{4}^{1/2}W)H_{\alpha }(W))
&=&\E(\Phi _{\delta _{n}}(I_{2}^{1/2}W^{\prime })h_{2}(W_{i}^{\prime
}))\times \E(\Phi _{\delta _{n}}(I_{2}^{1/2}W^{\prime \prime
})h_{2}(W_{j-2}^{\prime \prime })) \\
&\rightarrow & {}{\frac{1}{3}}\E(\left\vert B_{2}\right\vert \delta
_{0}(B_{1})h_{2}(B_{i}))\times \E(\left\vert B_{2}\right\vert \delta
_{0}(B_{1})h_{2}(B_{j-2}))
\end{eqnarray*}%
where $B=(B_{1},B_{2})$ is standard normal. If $i=1$ then

\begin{equation*}
\E(\left\vert B_{2}\right\vert \delta _{0}(B_{1})h_{2}(B_{1}))=\E(\left\vert
B_{2}\right\vert )\E(\delta _{0}(B_{1})(B_{1}^{2}-1))=-\frac{2}{\sqrt{2\pi }}%
\times \frac{1}{\sqrt{2\pi }}=-\frac{1}{\pi }
\end{equation*}%
and if $i=2$ then 
\begin{equation*}
\E(\left\vert B_{2}\right\vert \delta _{0}(B_{1})h_{2}(B_{2}))=\E(\left\vert
B_{2}\right\vert (B_{2}^{2}-1))\E(\delta _{0}(B_{1}))={}{\frac{1}{\pi}}.
\end{equation*}%
{}{So, discussing according to the possible values of $i,j$, we may define%
\begin{equation*}
\rho _{i,j}= \frac{1}{\pi^2} (-1)^{i+j}
\end{equation*}%
}
and we finally obtain, for $\alpha =(i,i,j,j)$%
\begin{equation*}
\lim_{n}\E(\partial _{1}\partial _{3}\Psi _{\delta
_{n}}(I_{4}^{1/2}W)H_{\alpha }(W))={}{\frac{1}{3}\rho _{i,j}}
\end{equation*}%
and 
{}{
\begin{eqnarray*}
\lim_{n}\frac{1}{n^{2}}\int_{D_{n,\varepsilon }}\gamma _{n}^{\prime
}(t,s)dsdt &=&6\times \frac{1}{2}\times \frac{1}{24}\sum_{i,j=1}^{2}\frac{1}{3}\rho _{i,j}\times \pi
^{2}U((i,i,j,j)) +O(\epsilon)\\
&=&\frac{1}{216}\sum_{i,j=1}^{2}\frac{(-3)^{i+j}}{4(1+2(i+j-2))}\times y_*+O(\epsilon)\\
&=&\frac{1}{120}\times y_*+O(\epsilon).
\end{eqnarray*}%
}

\textbf{Step 4}. We estimate $r_{n}(t,s)$ and $r_{n}(t).$ Since $L_{0}(\Phi
_{\delta _{n}})=1$ we have $\left\vert r_{n}(nt,ns)\right\vert \leq
\left\Vert \Sigma _{n}(nt,ns)-I_{4}\right\Vert .$ Let us compute $\Sigma
_{n}^{i,j}(nt,ns).$ By direct computations on has $\Sigma
_{n}^{1,1}(nt,ns)=\Sigma _{n}^{3,3}(nt,ns)=1$ and $\Sigma
_{n}^{1,2}(nt,ns)=\Sigma _{n}^{3,4}(nt,ns)=0.$ Moreover 
\begin{equation*}
\Sigma _{n}^{2,2}(nt,ns)=\frac{1}{n}\sum_{k=1}^{n}\E(Z_{n,k}^{2}(nt,Y))=\frac{%
1}{n}\sum_{k=1}^{n}\frac{k^{2}}{n^{2}}\rightarrow \int_{0}^{1}x^{2}dx=\frac{1%
}{3}=I_{4}^{2}.
\end{equation*}%
The same is true for $\Sigma _{n}^{4,4}(nt,ns).$ We look now to $\Sigma
_{n}^{i,j}(nt,ns)$ with $i\in \{1,2\}$ and $j\in \{3,4\}.$ Say for example
that $i=1$ and $j=4.$ Then we compute%
\begin{eqnarray*}
\E(Z_{n,k}^{1}(nt,Y)Z_{n,k}^{2}(ns,Y)) &=&\frac{k}{n}\E((\cos
(kt)Y_{k}^{1}+\sin (kt)Y_{k}^{2})(-\sin (ks)Y_{k}^{1}+\cos (ks)Y_{k}^{2})) \\
&=&\frac{k}{n}(\cos (ks)\sin (kt)-\cos (kt)\sin (ks))=\frac{k}{n}\sin
(k(t-s)).
\end{eqnarray*}%
Then, by using the ergodic lemma, if$\frac{t-s}{\pi }$ is irrational we get 
\begin{equation*}
\Sigma _{n}^{1,4}(nt,ns)=\frac{1}{n}\sum_{k=1}^{n}\frac{k}{n}\sin
(k(t-s))\rightarrow \frac{1}{4\pi }\int_{0}^{2\pi }\sin (u)du=0.
\end{equation*}%
The same result is obtained in the other cases. We conclude that $%
\lim_{n}r_{n}(nt,ns)=0$ $dt,ds$ almost surely. Since $\left\vert
r_{n}(nt,ns)\right\vert \leq 1,$ we may use Lebesgue's convergence theorem
and we obtain 
\begin{equation*}
\frac{1}{n^{2}}\int_{D_{\nu ,\varepsilon }}\left\vert r_{n}(t,s)\right\vert
dsdt=\int_{[0,\pi ]^{2}}1_{D_{n,\varepsilon }}(nt,ns)\left\vert
r_{n}(t,s)\right\vert dsdt\rightarrow 0.
\end{equation*}%
For $r_{n}(t)$ the same conclusion is (trivially) true.

\textbf{Step 5}. Estimate of $R_{n}(t,s,\Psi _{\delta _{n}}).$ Notice that $%
L_{0}(\Psi _{\delta _{n}})=1$ and $L_{q}(\Psi _{\delta _{n}})=\delta
_{n}^{-2}=n^{10}$ (see (\ref{N5b}))$.$ Then by (\ref{N7e}) 
\begin{equation*}
\left\vert R_{n}(t,s,\Psi _{\delta _{n}})\right\vert \leq C(1+n^{3/2}\times
n^{10}\times e^{-cn})
\end{equation*}%
with $C$ a constant which depends on $r,\eta $ from (\ref{D}) on $M_{p}(Y)$
from (\ref{M}) and on the lower eigenvalue $\varepsilon _{\ast }$ defined in
(\ref{ND}) for the covariance matrix $\Sigma _{n}(t,s).$ We have proved in (%
\ref{App5}) that this lower eigenvalue is lower bounded uniformly with
respect to $n$ so we conclude that the constant $C$ in the above inequality
does not depend on $n.$ Consequently 
\begin{equation*}
\sup_{n}\sup_{(t,s)\in D_{n,\varepsilon }}\left\vert R_{n}(t,s,\Psi _{\delta
_{n}})\right\vert \leq C<\infty
\end{equation*}%
and then 
\begin{equation*}
\frac{1}{n^{2}}\int_{D_{n,\varepsilon }}\frac{1}{\sqrt{n}}\left\vert
R_{n}(t,s,\Psi _{\delta _{n}})\right\vert dsdt\rightarrow 0.
\end{equation*}%
Similar estimates hold for $R_{n}(t,\Phi _{\delta _{n}}).$ Since $W$ is
standard normal, direct computations show that 
\begin{equation*}
\left\vert \E(\partial _{1}\Phi _{\delta _{n}}(I_{2}^{1/2}W)\Gamma
_{n,2}(I_{2}^{-1/2}Z_{n}(t,Y),W)))\right\vert \leq C
\end{equation*}%
and so 
\begin{equation*}
\frac{1}{n^{3}}\int_{D_{n,\varepsilon }}\left\vert \E(\partial _{1}\Phi
_{\delta _{n}}(I_{2}^{1/2}W)\Gamma
_{n,2}(I_{2}^{-1/2}Z_{n}(t,Y),W)))\E(\partial _{1}\Phi _{\delta
_{n}}(I_{2}^{1/2}W)\Gamma _{n,2}(I_{2}^{-1/2}Z_{n}(s,Y),W)))\right\vert
dsdt\rightarrow 0
\end{equation*}%
So we have proved that 
\begin{equation*}
\frac{1}{n^{2}}\int_{D_{n,\varepsilon }}\left\vert \overline{R}%
_{n}(t,s)\right\vert dsdt\rightarrow 0
\end{equation*}%
and the whole proof is completed.

$\square $

\appendix

\section{Ergodic lemma}
The following lemmas are based on the ergodic action of irrational rotations on the Torus.
\begin{lemma}
Set $\alpha$ a positive number such that $\frac{\alpha}{\pi}\in \mathbb{R}/%
\mathbb{Q}$, $f$ a $2\pi$--periodic function and $q\ge 1$ a positive
integer. One gets

\begin{equation}
\lim_{n\rightarrow \infty }\frac{1}{n}\sum_{k=1}^{n}f(k\alpha )=\frac{1}{%
2\pi }\int_{0}^{2\pi }f(t)dt,  \label{er1}
\end{equation}

\begin{equation}
\lim_{n\rightarrow \infty }\frac{1}{n}\sum_{k=1}^{n}\frac{k^{q}}{n^{q}}%
f(k\alpha )=\frac{1}{(q+1)2\pi }\int_{0}^{2\pi }f(t)dt  \label{e5}
\end{equation}
\end{lemma}

\textbf{Proof.} Let us denote by $\mathcal{C}_{2\pi }^{0}(\mathbb{R})$\ the
space of continuous $2\pi $ periodic functions and let introduce 
\begin{equation*}
\mathcal{H}_{0}=\left\{ \phi \in \mathcal{C}_{2\pi }^{0}(\mathbb{R})\,\,%
\Big{|}\,\,\int_{0}^{2\pi }\phi (t)dt=0\right\} ,
\end{equation*}%
and 
\begin{equation*}
\mathcal{E}=\left\{ f(x)=\phi (x+\alpha )-\phi (x)\,\,\Big{|}\,\,\phi \in 
\mathcal{C}_{2\pi }^{0}(\mathbb{R})\right\} .
\end{equation*}

Let us first prove that $\mathcal{E}$ is dense in $\mathcal{H}_{0}$. We take 
$T$ a continuous linear form on $\mathcal{H}_{0}$ and we extend it to $%
\mathcal{C}_{2\pi }^{0}(\mathbb{R})$ by taking $T(\phi )=T(\phi -m(\phi ))$
with $m(\phi )=\int_{0}^{2\pi }\phi (t)dt.$ We have to prove that if $T$
vanishes on $\mathcal{E}\ $then $T=0$ (in virtue of the Hahn-Banach
Theorem, this implies that $\mathcal{E}$ is dense in $\mathcal{H}_{0})$. The
Riesz Theorem ensures us that there exists a finite measure $\mu $ on $%
\mathbb{R}/2\pi \mathbb{Z}$ such that

\begin{equation*}
\forall \phi \in \mathcal{C}^0_{2\pi}(\mathbb{R}),\,\,T(\phi)=\int_0^{2\pi}
\phi(x) d\mu(x).
\end{equation*}

Since $Tf=0$ for every $f\in\mathcal{E}$, for any integer $n\ge 1$ one has

\begin{equation*}
\int_{0}^{2\pi} \phi(x +n \alpha) d\mu(x) =\int_0^{2\pi} \phi(x) d\mu(x),
\end{equation*}

and since the sequence $n\alpha$ is dense modulo $2\pi$ one deduces that for
any $y\in\mathbb{R}$:

\begin{equation*}
\int_{0}^{2\pi} \phi(x +y) d\mu(x) =\int_0^{2\pi} \phi(x) d\mu(x).
\end{equation*}
As a result, $\mu$ is invariant under translations and necessarily it is
Lebesgue measure up to a multiplicative constant. Hence, we get that $T=0$
over $\mathcal{H}_0$ and that $\mathcal{E}$ is dense for the uniform
topology. Finally, this preliminary consideration enables us to consider
that $f(x)=\phi(x+\alpha)-\phi(x)$ in the statements (\ref{er1}) and (\ref%
{e5}). Then, the conclusion is immediate since an Abel transforms gives us

\begin{equation*}
\Big|\frac{1}{n}\sum_{k=1} \frac{k^q}{n^q} \left(\varphi((k+1)\alpha)-%
\varphi(k\alpha)\right)\Big|\le 2 \|\phi\|_\infty\frac{1}{n} \sum_{k=1}^n
\left(\frac{(k+1)^q-k^q}{n^q}\right) \xrightarrow[n\to\infty]~0.
\end{equation*}
$\square$

In the following $C_{n}(k,t)$ is the matrix introduced in (\ref{N1}).

\begin{lemma}
\label{A}For every $i,j,l,l^{\prime }\in \{1,2\}$ and every $t,s$ such that $%
t,s,t+s,t-s$ are irrational one has 
\begin{equation}
\lim_{n}\frac{1}{n}\sum_{k=1}^{n}C_{n}^{i,l}(k,nt)^{2}C_{n}^{j,l^{\prime
}}(k,ns)^{2}=\frac{1}{4(1+2(i+j-2))}  \label{er9}
\end{equation}
\end{lemma}

\textbf{Proof.} We treat just two examples: take $i=1,j=2,l=1,l^{\prime }=2.$
Then%
\begin{eqnarray*}
C_{n}^{i,l}(k,nt)^{2}C_{n}^{j,l^{\prime }}(k,ns)^{2} &=&(\cos kt\times \frac{
k}{n}\cos ks)^{2}=\frac{1}{4}\times \frac{k^{2}}{n^{2}}(\cos (k(t+s))+\cos
(k(t-s)))^{2} \\
&=&\frac{1}{4}\times \frac{k^{2}}{n^{2}}(\cos ^{2}(k(t+s))+\cos
^{2}(k(t-s))+2\cos (k(t+s))\cos (k(t-s))) \\
&=&\frac{1}{4}\times \frac{k^{2}}{n^{2}}(\cos ^{2}(k(t+s))+\cos
^{2}(k(t-s))+\cos (2kt)-\cos (2ks))).
\end{eqnarray*}%
Then, the ergodic lemma (with $q=2)$ gives 
\begin{equation*}
\lim_{n}\frac{1}{n}\sum_{k=1}^{n}C_{n}^{i,l}(k,nt)^{2}C_{n}^{j,l^{\prime
}}(k,ns)^{2}=2\times \frac{1}{4}\times \frac{1}{2\pi \times 3}\int_{0}^{2\pi
}(\cos ^{2}(u)+\cos (u))du=\frac{1}{12}.
\end{equation*}

Take now $i=2,j=2,l=1,l^{\prime }=2.$ Then%
\begin{eqnarray*}
C_{n}^{i,l}(k,nt)^{2}C_{n}^{j,l^{\prime }}(k,ns)^{2} &=&(\frac{k}{n}\sin
kt\times \frac{k}{n}\cos ks)^{2}=\frac{1}{4}\times \frac{k^{4}}{n^{4}}(\sin
(k(t+s))+\sin (k(t-s))^{2} \\
&=&\frac{1}{4}\times \frac{k^{4}}{n^{4}}(\sin ^{2}(k(t+s))+\sin
^{2}(k(t-s))+\cos (2kt)+\cos (2ks))
\end{eqnarray*}%
Then, the ergodic lemma (with $q=4)$ gives 
\begin{equation*}
\lim_{n}\frac{1}{n}\sum_{k=1}^{n}C_{n}^{i,l}(k,nt)^{2}C_{n}^{j,l^{\prime
}}(k,ns)^{2}=\frac{1}{4}\times \frac{1}{2\pi \times 5}\times 2\int_{0}^{2\pi
}(\sin ^{2}(u)+\cos (u))du=\frac{1}{20}.
\end{equation*}%
$\square $

\bigskip

\begin{lemma}
\label{B}For every $j,i,l\in \{1,2\}$ and every $t,s$ such that $t,s,t+s,t-s$
are irrational one has 
\begin{eqnarray}
&&\lim_{n}\frac{1}{n}
\sum_{k=1}^{n}C_{n}^{i,1}(k,nt)C_{n}^{i,2}(k,nt)C_{n}^{j,l}(k,ns)^{2}
\label{er10} \\
&=&\lim_{n}\frac{1}{n}
\sum_{k=1}^{n}C_{n}^{i,1}(k,nt)C_{n}^{i,2}(k,nt)C_{n}^{j,1}(k,ns)C_{n}^{j,2}(k,ns)=0.
\notag
\end{eqnarray}
\end{lemma}

\textbf{Proof.} All the computations are analogous so we treat just an
example: $l=i=j=1.$ So we have%
\begin{eqnarray*}
C_{n}^{1,1}(k,nt)C_{n}^{1,2}(k,nt)C_{n}^{1,1}(k,ns)^{2} &=&\cos (kt)\sin
(kt)\cos ^{2}(ks) \\
&=&\frac{1}{2}\sin (2kt)\cos ^{2}(ks) \\
&=&\frac{1}{4}(\sin (k(2t+s))+\sin (k(2t-s)))\cos (ks) \\
&=&\frac{1}{8}(\sin (k(2t+2s))+2\sin (2kt)+\sin (k(2t-2s)))
\end{eqnarray*}%
and using the ergodic lemma with $q=0$ we get%
\begin{equation*}
\lim_{n}\frac{1}{n}%
\sum_{k=1}^{n}C_{n}^{1,1}(k,nt)C_{n}^{1,2}(k,nt)C_{n}^{1,1}(k,ns)^{2}=\frac{1%
}{8}\times \frac{1}{2\pi }\times 4\int_{0}^{2\pi }\sin (u)du=0.
\end{equation*}%
$\square $

\begin{lemma}
\label{C}For every $i_{1},i_{2},i_{3},l_{1},l_{2},l_{3}\in \{1,2\}$ and
every $t,s$ such that $t,s,t+s,t-s$ are irrational one has 
\begin{equation}
\lim_{n}\frac{1}{n}
\sum_{k=1}^{n}C_{n}^{i_{1},l_{1}}(k,nt)C_{n}^{i_{2},l_{2}}(k,nt)C_{n}^{i_{3},l_{3}}(k,ns)=0.
\label{er11}
\end{equation}
\end{lemma}

\textbf{Proof.} The poof is similar in all cases so we treat just an
example: $i_{1}=1,i_{2}=2,i_{3}=2,l_{1}=l_{2}=l_{3}=1.$ Then we deal with%
\begin{eqnarray*}
\cos (kt)\times \frac{k}{n}\sin (kt)\times \frac{k}{n}\sin (ks) &=&\frac{
k^{2}}{n^{2}}\times \frac{1}{2}\sin (2kt)\sin (ks) \\
&=&\frac{k^{2}}{n^{2}}\times \frac{1}{4}(\cos (k(2t+s))-\cos (k(2t-s)))
\end{eqnarray*}%
And using the ergodic lemma with $q=2$ we get%
\begin{eqnarray*}
&&\lim_{n}\frac{1}{n}\sum_{k=1}^{n}\frac{k^{2}}{n^{2}}\times \frac{1}{4}
(\cos (k(2t+s))-\cos (k(2t-s))) \\
&=&\frac{1}{24\pi }(\int_{0}^{2\pi }\cos (u)du-\int_{0}^{2\pi }\cos (u)du)=0.
\end{eqnarray*}%
$\square $

\section{Estimates of some trigonometric sums}

For $n\in N,$ $i=0,1,2$ and $b\in R_{+}\backprime \{2\pi p;p\in N\}$ we put 
\begin{equation*}
S_{b,i}(c)=\frac{1}{n}\sum_{k=1}^{n}\frac{k^{i}}{n^{i}}\cos (bk).
\end{equation*}%
We also denote 
\begin{equation}
\overline{b}=\inf_{p\in N}\frac{\left\vert 2\pi p-b\right\vert }{p\vee 1}.
\label{c1}
\end{equation}%
The aim of this section is to prove the following lemma:

\begin{lemma}
There exists an universal constant $C\geq 1$ such that for every $n\in N$ $%
i=0,1,2$ and $\in R_{+}\backprime \{2\pi p;p\in N\}$%
\begin{equation}
\left\vert S_{b,i}(c)\right\vert \leq \frac{C}{n\overline{b}}.  \label{c1'}
\end{equation}
\end{lemma}

The first step is the following abstract estimate:

\begin{lemma}
\textbf{A}. Let $f\in L^{2}(0,1)$ and let $\phi (x)=\sum_{k=0}^{\infty
}f(x-k)1_{[k,k+1)}(x).$ There exists an universal constant such that for
every $k<n$%
\begin{equation}
\left\vert \int_{k}^{n}\phi (x)\cos (bx)\right\vert \leq \frac{C}{\overline{b%
}}\left\Vert f\right\Vert _{2}.  \label{c2}
\end{equation}%
\textbf{B}. Moreover there exists an universal constant $C$ such that, for $%
i=0,1,2$%
\begin{equation}
\left\vert \int_{0}^{n}\frac{x^{i}}{n^{i}}\phi (x)\cos (bx)dx\right\vert
\leq \frac{C}{\overline{b}}\left\Vert f\right\Vert _{2}  \label{c3}
\end{equation}%
and in particular, taking $f=1,$ 
\begin{equation}
\int_{0}^{n}\frac{x^{i}}{n^{i}}\cos (bx)dx\leq \frac{C}{\overline{b}}
\label{c3'}
\end{equation}%
The same estimates hold if we replace $\cos $ by $\sin .$
\end{lemma}

\textbf{Proof of A.} We denote $\alpha _{0}=\int_{0}^{1}f(x)dx$ and 
\begin{equation*}
\alpha _{p}=\int_{0}^{1}f(x)\cos (2\pi px)dx,\quad \beta
_{p}=\int_{0}^{1}f(x)\sin (2\pi px)dx.
\end{equation*}%
Then using the development in Fourier series of $\phi $ we obtain 
\begin{eqnarray*}
\int_{k}^{n}\phi (x)\cos (bx)dx &=&\alpha _{0}\int_{k}^{n}\cos (bx)dx \\
&&+\sum_{p=1}^{\infty }\alpha _{p}\int_{k}^{n}\cos (2\pi px)\cos
(bx)dx+\beta _{p}\int_{k}^{n}\sin (2\pi px)\cos (bx)dx.
\end{eqnarray*}%
We write%
\begin{equation*}
\cos (2\pi px)\cos (bx)=\frac{1}{2}(\cos ((2\pi p+b)x)+\cos ((2\pi p-b)x).
\end{equation*}%
and we use a similar decomposition for $\sin (2\pi px)\cos (bx).$

Notice that for every $\theta >0$ one has 
\begin{equation*}
\left\vert \int_{k}^{n}\cos (\theta x)dx\right\vert \leq \frac{2\pi }{\theta 
}\quad and\quad \left\vert \int_{k}^{n}\sin (\theta x)dx\right\vert \leq 
\frac{2\pi }{\theta }.
\end{equation*}%
Using these inequalities we obtain 
\begin{eqnarray*}
\left\vert \int_{k}^{n}\phi (x)\cos (bx)\right\vert &=&\frac{2\pi }{b}%
\left\vert \alpha _{0}\right\vert +\sum_{p=1}^{\infty }(\left\vert \alpha
_{p}\right\vert +\left\vert \beta _{p}\right\vert (\frac{2\pi }{2\pi p+b}+%
\frac{2\pi }{\left\vert 2\pi p-b\right\vert }) \\
&\leq &\frac{2\pi }{b}\left\vert \alpha _{0}\right\vert +\frac{4\pi }{%
\overline{b}}\sum_{p=1}^{\infty }(\left\vert \alpha _{p}\right\vert
+\left\vert \beta _{p}\right\vert )\frac{1}{p} \\
&\leq &\frac{2\pi }{b}\left\vert \alpha _{0}\right\vert +\frac{C}{\overline{b%
}}(\sum_{p=1}^{\infty }(\left\vert \alpha _{p}\right\vert +\left\vert \beta
_{p}\right\vert )^{2})^{1/2}\leq \frac{C}{\overline{b}}\left\Vert
f\right\Vert _{2}.
\end{eqnarray*}

\textbf{B}. We just treat the case $i=1$ (the other ones are similar). We
write%
\begin{equation*}
\frac{1}{n}\int_{0}^{n}x\phi (x)\cos (bx)dx=\frac{1}{n}\int_{0}^{n}\psi
(x)\cos (bx)dx+\frac{1}{n}\sum_{k=1}^{n}k\int_{k}^{k+1}\phi (x)\cos (bx)dx
\end{equation*}%
with $\psi $ associated to $g($ $x)=xf(x).$ Using (\ref{c2}) (notice that $%
\left\Vert g\right\Vert _{2}\leq \left\Vert f\right\Vert _{2})$ 
\begin{equation*}
\frac{1}{n}\left\vert \int_{0}^{n}\psi (x)\cos (bx)dx\right\vert \leq \frac{C%
}{n}\times \frac{1}{\overline{b}}\left\Vert f\right\Vert _{2}.
\end{equation*}%
Moreover 
\begin{equation*}
\frac{1}{n}\sum_{k=1}^{n}k\int_{k}^{k+1}\phi (x)\cos (bx)dx=\frac{1}{n}%
\sum_{k=1}^{n}\int_{k}^{n}\phi (x)\cos (bx)dx
\end{equation*}%
so that, by (\ref{c2}) we upper bound the above term by 
\begin{equation*}
\frac{1}{n}\sum_{k=1}^{n}\left\vert \int_{k}^{n}\phi (x)\cos
(bx)dx\right\vert \leq \frac{C}{n}\times n\times \frac{1}{\overline{b}}%
\left\Vert f\right\Vert _{2}=\frac{C}{\overline{b}}\left\Vert f\right\Vert
_{2}.
\end{equation*}%
And (\ref{c3'}) is a particular case of (\ref{c3}) with $f=1$ (so that $\phi
=1).$ $\square $

We recall that 
\begin{equation*}
S_{b,i}(c)=\frac{1}{n}\sum_{k=1}^{n}\frac{k^{i}}{n^{i}}\cos (\frac{ak}{n}%
),\quad with\quad a=nb
\end{equation*}%
We also denote%
\begin{eqnarray}
\phi (x) &=&\sum_{k=0}^{\infty }(k+1-x)^{2}1_{\{k\leq x<k+1\}}\quad and\quad
J_{b,i}(c)=\frac{1}{n^{1+i}}\int_{0}^{n}x^{i}\phi (x)\cos (bx)dx  \label{r0}
\\
I_{b,i}(c) &=&\frac{1}{n^{1+i}}\int_{0}^{n}x^{i}\cos
(bx)dx=\int_{0}^{1}x^{i}\cos (ax)dx.
\end{eqnarray}%
We also define $S_{b,i}(s)$ and $J_{b,i}(s)$ by replacing the $\cos $ by $%
\sin .$

We will prove the following estimates:

\begin{lemma}
Let $a=bn$ with $0<b.$ There exists an universal constant $C$ such that, for 
$i=0,1,2$%
\begin{equation}
\left\vert S_{b,i}(c)\right\vert =\frac{1}{1+b^{2}/4}(\frac{b^{3}}{4}%
J_{b,i}(s)+\frac{b^{2}}{2}J_{b,i}(c)+I_{b,i}(c)-\frac{b}{2}%
I_{b,i}(s))+\varepsilon _{n}  \label{r1}
\end{equation}%
with $\left\vert \varepsilon _{n}\right\vert \leq C/n.$
\end{lemma}

\textbf{Proof.} Let us prove (\ref{r1}) for $i=1$ (the proof is analogous
for $i=0$ and $i=2)$. We write%
\begin{equation*}
S_{b,1}(c)=I_{b,1}(c)-\sum_{k=1}^{n}\int_{k/n}^{(k+1)/n}(x\cos (ax)-\frac{k}{%
n}\cos (\frac{ak}{n}))dx.
\end{equation*}%
Moreover 
\begin{equation*}
\int_{k/n}^{(k+1)/n}(x\cos (ax)-\frac{k}{n}\cos (\frac{ak}{n}))dx=\frac{k}{n}%
\int_{k/n}^{(k+1)/n}(\cos (ax)-\cos (\frac{ak}{n}))dx+\delta _{n,k}
\end{equation*}%
with $\left\vert \delta _{n,k}\right\vert \leq 1/n^{2}$ so that $%
\sum_{k=1}^{n}\delta _{n,k}=\varepsilon _{n},$ with $\left\vert \varepsilon
_{n}\right\vert \leq C/n.$ We write now (recall that $a=nb)$%
\begin{eqnarray*}
\frac{k}{n}\int_{k/n}^{(k+1)/n}(\cos (ax)-\cos (\frac{ak}{n}))dx &=&-\frac{ak%
}{n}\int_{k/n}^{(k+1)/n}dx\int_{k/n}^{x}dy\sin ay \\
&=&-\frac{ak}{2n^{3}}\sin \frac{ak}{n}-\frac{ak}{n}\int_{k/n}^{(k+1)/n}dx%
\int_{k/n}^{x}dy(\sin ay-\sin \frac{ak}{n}) \\
&=&-\frac{bk}{2n^{2}}\sin \frac{ak}{n}-\frac{a^{2}k}{n}%
\int_{k/n}^{(k+1)/n}dx\int_{k/n}^{x}dy\int_{k/n}^{y}dz\cos az \\
&=&-\frac{bk}{2n^{2}}\sin \frac{ak}{n}-\frac{a^{2}k}{2n}%
\int_{k/n}^{(k+1)/n}(k+1-z)^{2}\cos (az)dz
\end{eqnarray*}%
Summing over $k$ this gives (with $\varepsilon _{n}$ of order $\frac{1}{n}$
and which changes from a line to another)%
\begin{eqnarray*}
S_{b,1}(c) &=&I_{b,1}(c)+\varepsilon _{n}+\frac{b}{2}S_{b,1}(s)+\frac{a^{2}}{%
2}\sum_{k=1}^{n}\frac{k}{n}\int_{k/n}^{(k+1)/n}\cos (ax)(\frac{k+1}{n}-x)^{2}
\\
&=&I_{b,1}(c)+\varepsilon _{n}+\frac{b}{2}S_{b,1}(s)+\frac{a^{2}}{2}%
\sum_{k=1}^{n}\int_{k/n}^{(k+1)/n}x\cos (ax)(\frac{k+1}{n}-x)^{2} \\
&=&I_{b,1}(c)+\varepsilon _{n}+\frac{b}{2}S_{b,1}(s)+\frac{b^{2}}{2}\times 
\frac{1}{n^{2}}\int_{0}^{n}x\phi (x)\cos (bx)dx \\
&=&I_{b,1}(c)+\varepsilon _{n}+\frac{b}{2}S_{b,1}(s)+\frac{b^{2}}{2}%
J_{b,1}(c).
\end{eqnarray*}%
The same computations give 
\begin{equation*}
S_{n,1}(s)=I_{b,1}(s)+\varepsilon _{n}-\frac{b}{2}S_{n,1}(c)+\frac{b^{2}}{2}%
J_{b,1}(s).
\end{equation*}%
We insert this in the previous estimate and we get 
\begin{equation*}
S_{b,1}(c)=\varepsilon _{n}+I_{1}(c)-\frac{b}{2}I_{b,1}(s)-\frac{b^{2}}{4}%
S_{b,1}(c)+\frac{b^{3}}{4}J_{b,1}(s)+\frac{b^{2}}{2}J_{b,1}(c)
\end{equation*}%
and we are done. $\square $

\bigskip

\textbf{Proof of (\ref{c1'}).} By (\ref{c3}) and (\ref{c3'}) 
\begin{equation*}
\left\vert J_{b,i}(s)\right\vert +\left\vert J_{b,i}(c)\right\vert
+\left\vert I_{b,i}(s)\right\vert +\left\vert I_{b,i}(c)\right\vert \leq 
\frac{C}{n\overline{b}}
\end{equation*}%
so \ (\ref{c1'}) follows. $\square $

\section{Non degeneracy}

In this section we discuss the non degeneracy of the matrix $\Sigma
_{n}(t,s) $ which is the covariance matrix of $S_{n}(t,s).$ Direct
computations show that: 
\begin{eqnarray}
\Sigma _{n}^{1,1}(t,s) &=&\Sigma _{n}^{3,,3}(t,s)=1,\quad \Sigma
_{n}^{2,2}(t,s)=\Sigma _{n}^{4,4}(t,s)=\frac{1}{n}\sum_{k=1}^{n}\frac{k^{2}}{%
n^{2}},\quad  \label{App3} \\
\Sigma _{n}^{1,3}(t,s) &=&\Sigma _{n}^{3,1}(t,s)=\frac{1}{n}%
\sum_{k=1}^{n}\cos \frac{k(t-s)}{n}\quad \Sigma _{n}^{2,4}(t,s)=\Sigma
_{n}^{4,2}(t,s)=\frac{1}{n}\sum_{k=1}^{n}\frac{k^{2}}{n^{2}}\cos \frac{k(t-s)%
}{n}  \notag \\
\Sigma _{n}^{1,4}(t,s) &=&\Sigma _{n}^{4,1}(t,s)=-\Sigma
_{n}^{2,3}(t,s)=-\Sigma _{n}^{3,2}(t,s)=\frac{1}{n}\sum_{k=1}^{n}\frac{k}{n}%
\cos \frac{k(t-s)}{n}  \notag \\
\Sigma _{n}^{1,2}(t,s) &=&\Sigma _{n}^{2,1}(t,s)=\Sigma
_{n}^{3,4}(t,s)=\Sigma _{n}^{4,3}(t,s)=0.  \notag
\end{eqnarray}%
We define $\Sigma (t,s)$ just by passing to the limit (for fixed $t$ and $%
s): $%
\begin{eqnarray*}
\Sigma ^{1,1}(t,s) &=&\Sigma ^{3,,3}(t,s)=1,\quad \Sigma ^{2,2}(t,s)=\Sigma
^{4,4}(t,s)=\frac{1}{3} \\
\Sigma _{n}^{1,3}(t,s) &=&\Sigma ^{3,1}(t,s)=\int_{0}^{1}\cos
((t-s)x)dx\quad \Sigma ^{2,4}(t,s)=\Sigma ^{4,2}(t,s)=\int_{0}^{1}x^{2}\cos
((t-s)x)dx \\
\Sigma ^{1,4}(t,s) &=&\Sigma ^{4,1}(t,s)=-\Sigma ^{2,3}(t,s)=-\Sigma
^{3,2}(t,s)=\int_{0}^{1}x\cos ((t-s)x)dx \\
\Sigma ^{1,2}(t,s) &=&\Sigma ^{2,1}(t,s)=\Sigma ^{3,4}(t,s)=\Sigma
^{4,3}(t,s)=0.
\end{eqnarray*}%
Then it is easy to check that there exists an universal constant $C\geq 1$
such that for every $i,j=1,...,4$ and every $0<s<t$%
\begin{equation}
\sup_{\left\vert t-s\right\vert \leq n^{\rho }}\left\vert \Sigma
_{n}^{i,j}(t,s)-\Sigma (t,s)\right\vert \leq \frac{C(t-s)}{n}.  \label{App1}
\end{equation}%
Notice however that, if $t-s\sim n$ the above inequality says nothing. So
our strategy will be the following: we consider a first case, when $t-s\leq 
\sqrt{n}$ and then we use that non degeneracy of $\Sigma (t,s)$ (which we
prove in the following lemma) in order to obtain the non degeneracy of $%
\Sigma _{n}(t,s).$ And in the case $\sqrt{n}\leq t-s\leq n\pi $ we use the
estimates from the previous section in order to obtain directly the non
degeneracy of $\Sigma _{n}(t,s)$.

\begin{lemma}
\label{NonDeg}For every $\varepsilon >0$ there exists some $\lambda
(\varepsilon )>0$ such that for every $t$ and $s$ such that $\left\vert
t-s\right\vert >\varepsilon $ one has 
\begin{equation}
\det \Sigma (t,s)\geq \lambda (\varepsilon ).  \label{App4}
\end{equation}
\end{lemma}

\textbf{Proof.} Using integration by parts%
\begin{equation*}
\Sigma ^{2,4}(t,s)=\frac{1}{t-s}(\sin (t-s)-2\int_{0}^{1}x\sin ((t-s)x)dx
\end{equation*}%
so $\Sigma ^{2,4}(t,s)\rightarrow 0$ as $t-s\rightarrow \infty .$ The same
is true for every $\Sigma ^{i,j}(t,s)$ with $i\neq j.$ It follows that $%
\lim_{t-s\rightarrow \infty }\det \Sigma (t,s)=\frac{1}{9}.$

We will prove that for every $s<t$ and $\xi \neq 0$ one has%
\begin{equation}
\left\langle \Sigma (t,s)\xi ,\xi \right\rangle >0.  \label{App4'}
\end{equation}%
This implies that $\det \Sigma (t,s)>0$ and, since $t-s\rightarrow \det
\Sigma (t,s)$ is a continuous function, it has a strictly positive infimum $%
\lambda (\varepsilon )>0,$ so we obtain (\ref{App4}).

We will use the notation 
\begin{equation*}
C_{t}(x)=\cos (tx)\quad S_{t}(x)=\sin (tx).
\end{equation*}

\textbf{Step 1}. To begin we compute $\left\langle \Sigma (t,s)\xi ,\xi
\right\rangle .$ We have 
\begin{equation*}
\left\langle \Sigma _{n}(t,s)\xi ,\xi \right\rangle =\frac{1}{n}%
\sum_{i,j=1}^{4}\sum_{k=1}^{n}\E(\xi _{i}Z_{n,k}^{i}(t,s)\xi
_{j}Z_{n,k}^{j}(t,s))=\frac{1}{n}\sum_{k=1}^{n}\E((\sum_{i=1}^{4}\xi
_{i}Z_{n,k}^{i}(t,s))^{2}).
\end{equation*}%
We compute now 
\begin{eqnarray*}
\sum_{i=1}^{4}\xi _{i}Z_{n,k}^{i}(t,s) &=&\xi _{1}Z_{n,k}^{1}(t)+\xi
_{2}Z_{n,k}^{2}(t)+\xi _{3}Z_{n,k}^{1}(s)+\xi _{4}Z_{n,k}^{2}(s) \\
&=&\xi _{1}(Y_{k}^{1}C_{t}(\frac{k}{n})+Y_{k}^{2}S_{t}(\frac{k}{n}))+\xi _{2}%
\frac{k}{n}(-Y_{k}^{1}S_{t}(\frac{k}{n})+Y_{k}^{2}C_{t}(\frac{k}{n})) \\
&&+\xi _{3}(Y_{k}^{1}C_{s}(\frac{k}{n})+Y_{k}^{2}S_{s}(\frac{k}{n}))+\xi _{4}%
\frac{k}{n}(-Y_{k}^{1}S_{s}(\frac{k}{n})+Y_{k}^{2}C_{s}(\frac{k}{n})) \\
&=&Y_{k}^{1}\times (\xi _{1}C_{t}(\frac{k}{n})-\xi _{2}\frac{k}{n}S_{t}(%
\frac{k}{n})+\xi _{3}C_{s}(\frac{k}{n})-\xi _{4}\frac{k}{n}S_{s}(\frac{k}{n}%
)) \\
&&+Y_{k}^{2}\times (\xi _{1}S_{t}(\frac{k}{n})+\xi _{2}\frac{k}{n}C_{t}(%
\frac{k}{n})+\xi _{3}S_{s}(\frac{k}{n})+\xi _{4}\frac{k}{n}C_{s}(\frac{k}{n}%
))
\end{eqnarray*}%
so that, using the orthogonality of $Y_{k}^{1}$ and $Y_{k}^{2}$%
\begin{eqnarray*}
\E(\sum_{i=1}^{4}\xi _{i}Z_{n,k}^{i}(t,s))^{2} &=&(\xi _{1}C_{t}(\frac{k}{n}%
)-\xi _{2}\frac{k}{n}S_{t}(\frac{k}{n})+\xi _{3}C_{s}(\frac{k}{n})-\xi _{4}%
\frac{k}{n}S_{s}(\frac{k}{n}))^{2} \\
&&+(\xi _{1}S_{t}(\frac{k}{n})+\xi _{2}\frac{k}{n}C_{t}(\frac{k}{n})+\xi
_{3}S_{s}(\frac{k}{n})+\xi _{4}\frac{k}{n}C_{s}(\frac{k}{n}))^{2}.
\end{eqnarray*}%
By passing to the limit ($t,s$ and $\xi $ are fixed)%
\begin{equation*}
\left\langle \Sigma (t,s)\xi ,\xi \right\rangle
=\int_{0}^{1}(I_{1}^{2}(x)+I_{2}^{2}(x))dx
\end{equation*}%
with%
\begin{eqnarray*}
I_{1}(x) &=&\xi _{1}C_{t}(x)-\xi _{2}xS_{t}(x)+\xi _{3}C_{s}(x)-\xi
_{4}xS_{s}(x) \\
I_{2}(x) &=&\xi _{1}S_{t}(x)+\xi _{2}xC_{t}(x)+\xi _{3}S_{s}(x)+\xi
_{4}xC_{s}(x)
\end{eqnarray*}

\textbf{Step 2}. Suppose that for some $s<t$ and $\xi $ we have $%
\left\langle \Sigma (t,s)\xi ,\xi \right\rangle =0.$\ This implies that $%
I_{1}(x)=I_{2}(x)=0$ for every $x\in \lbrack 0,1].$ We will prove that $%
I_{1}(x)=0$ for every $x\in (0,1)$ then either $t=s$ or $\xi =0.$ Using the
developments in Taylor series for $\cos $ and $\sin $ we get%
\begin{eqnarray*}
C_{t}(x) &=&\cos (tx)=\sum_{n=0}^{\infty }\frac{(-1)^{n}}{(2n)!}x^{2n}\times
t^{2n} \\
xS_{t}(x) &=&x\sin (tx)=-\sum_{n=1}^{\infty }\frac{(-1)^{n}}{(2n)!}%
x^{2n}\times (2nt^{2n-1}).
\end{eqnarray*}%
We use this development in order to compute $I_{1}(x).$ This will be a
Taylor series and, if $I_{1}(x)=0$ for every $x\in (0,1)$ then it has null
coefficients. For $n=0$ we obtain 
\begin{equation*}
\xi _{1}=-\xi _{3}.
\end{equation*}%
We write then the equalities for $n=1,2,3:$%
\begin{equation*}
\xi _{1}(t^{2n}-s^{2n})=2n(\xi _{2}t^{2n-1}+\xi _{4}s^{2n-1}).
\end{equation*}%
Using the equation for $n=1$ we get the following equalities (corresponding
to $n=2$ respectively to $n=3)$ 
\begin{eqnarray*}
2(\xi _{2}t^{3}+\xi _{4}s^{3}) &=&(\xi _{2}t+\xi _{4}s)(t^{2}+s^{2}) \\
3(\xi _{2}t^{5}+\xi _{4}s^{5}) &=&(\xi _{2}t+\xi
_{4}s)(t^{4}+t^{2}s^{2}+s^{4}).
\end{eqnarray*}%
The first equation gives%
\begin{equation*}
\xi _{4}=\xi _{2}\times \frac{t}{s}
\end{equation*}%
and if we insert this in the second equation we obtain%
\begin{equation*}
3\xi _{2}(t^{4}+s^{4})=2\xi _{2}(t^{4}+t^{2}s^{2}+s^{4})
\end{equation*}%
which finally gives%
\begin{equation*}
\xi _{2}(t^{2}-s^{2})^{2}=0.
\end{equation*}%
$\square $

\begin{remark}
\label{cos}The above estimate of the covariance matrix is based on $%
Y_{k}^{1} $ only. This means that we obtain the invariance principle for
series of $\cos $ only (taking $Y_{k}^{2}=0)$ The reasoning is similar. But
we have to assume that $Y_{k}^{1}$ satisfies the Doeblin's condition on $R$
instead of $Y_{k}=(Y_{k}^{1},Y_{k}^{2})$ which is assumed in the present
paper to verify the Doeblin's condition in $R^{2}.$ And the correctors will
also change, so there is some work to do. We leave out this problem here.
\end{remark}

\begin{corollary}
\label{ND}Let $b_{\ast }<2\pi .$ For every $\varepsilon >0$ there exists $%
n(\varepsilon )$ such that for $n\geq n(\varepsilon )$ one has 
\begin{equation}
\inf_{\varepsilon <\left\vert t-s\right\vert \leq b_{\ast }n}\det \Sigma
_{n}(t,s)\geq \frac{1}{2}\lambda (\varepsilon )  \label{App5}
\end{equation}%
with $\lambda (\varepsilon )$ from (\ref{App4}).
\end{corollary}

\textbf{Proof.} Suppose first that $\varepsilon <t-s\leq n^{1/2}.$ Then 
\begin{equation*}
\det \Sigma _{n}(t,s)\geq \det \Sigma (t,s)-\left\vert \det \Sigma
_{n}(t,s)-\det \Sigma (t,s)\right\vert \geq \lambda (\varepsilon )-\frac{C}{%
n^{1/2}}\geq \frac{1}{2}\lambda (\varepsilon )
\end{equation*}%
for sufficiently large $n.$

We consider now the case $t-s>n^{1/2}.$ We will use (\ref{c1'}) with $b=%
\frac{t-s}{n}$ in order to prove that all the terms out of the diagonal are
very small, so the determinant will be close to the product of the terms of
the diagonal which is (almost) $\frac{1}{9}.$ We look to 
\begin{equation*}
\Sigma _{n}^{4,2}(t,s)=\frac{1}{n}\sum_{k=1}^{n}\frac{k^{2}}{n^{2}}\cos 
\frac{k(t-s)}{n}=S_{b,2}(c)
\end{equation*}%
Since $t-s\leq b_{\ast }n$ it follows that $b=\frac{t-s}{n}\leq b_{\ast
}<2\pi $ and this guarantees that $\overline{b}=\min \{b,2\pi -b_{\ast }\}.$
Since $nb=t-s\geq \sqrt{n},$ for sufficiently large $n$ we have $\overline{b}%
n\geq \sqrt{n}$ and so, by (\ref{c1'})%
\begin{equation*}
\left\vert \Sigma _{n}^{4,2}(t,s)\right\vert \leq \frac{C}{\sqrt{n}}%
\rightarrow 0.
\end{equation*}%
The same is true for the other terms out of the diagonal. $\square $

\end{document}